\documentclass{amsart}
\usepackage{amsmath,amssymb,amsthm,mathtools}
\usepackage[margin=1in]{geometry}
\usepackage{enumitem}
\usepackage{needspace}
\usepackage[colorlinks=true]{hyperref}
\usepackage{aliascnt}
\usepackage[capitalize,nameinlink]{cleveref}

\newtheorem{theorem}{Theorem}[section]
\newaliascnt{proposition}{theorem}
\newtheorem{proposition}[proposition]{Proposition}
\aliascntresetthe{proposition}
\newaliascnt{lemma}{theorem}
\newtheorem{lemma}[lemma]{Lemma}
\aliascntresetthe{lemma}
\newaliascnt{corollary}{theorem}
\newtheorem{corollary}[corollary]{Corollary}
\aliascntresetthe{corollary}
\theoremstyle{definition}
\newaliascnt{definition}{theorem}
\newtheorem{definition}[definition]{Definition}
\aliascntresetthe{definition}
\newaliascnt{assumption}{theorem}
\newtheorem{assumption}[assumption]{Assumption}
\aliascntresetthe{assumption}
\crefname{assumption}{Assumption}{Assumptions}
\Crefname{assumption}{Assumption}{Assumptions}
\theoremstyle{remark}
\newaliascnt{remark}{theorem}
\newtheorem{remark}[remark]{Remark}
\aliascntresetthe{remark}

\numberwithin{equation}{section}

\newcommand{\binomset}[2]{\binom{#1}{#2}}
\newcommand{\one}{\mathbf 1}
\newcommand{\ceff}{c_{\mathrm{eff}}}
\newcommand{\cL}{\mathcal L}
\newcommand{\cK}{\mathcal K}
\newcommand{\cH}{\mathcal H}
\newcommand{\cE}{\mathcal E}
\newcommand{\RR}{\mathbb R}

\title[Rigidity of the moving particle lemma]{Rigidity of the moving particle lemma for the interchange process}
\author{Joe P.\@ Chen}
\address{Department of Mathematics, Colgate University, Hamilton, NY 13346, USA}
\email{jpchen@colgate.edu}
\date{\today}
\keywords{Interchange process, exclusion process, moving particle lemma, octopus inequality, electric networks, representations of the symmetric group.}
\subjclass[2020]{
05C50, %Graphs and linear algebra (matrices, eigenvalues, etc.)
20C30, %Representations of finite symmetric groups
60J45, %Probabilistic potential theory
60K35. %Interacting random processes; statistical mechanics type models; percolation theory
}

\begin{document}

\begin{abstract}
Let $G=(V,c)$ be a finite connected weighted graph.
We establish the exact conditions for attaining equality in the moving particle lemma for the interchange process on $G$ subject to a unit voltage drop across two terminal vertices $x\neq y \in V$.  
First, equality can occur only in irreducible representations indexed by partitions with at most two rows.  
We then show that the equality phenomenon beyond the constant and linear sectors is already visible in the two-particle exclusion process. 
Such equality occurs in exactly two cases: either the harmonic $x$--$y$ voltage is $\frac12$ on every nonterminal
vertex, and the induced nonterminal graph $G[V\setminus\{x,y\}]$ is disconnected; or the harmonic voltage is constant on each connected component of $G[V\setminus \{x,y\}]$, and both voltage values $a$ and $1-a$, with $0\le a<\frac12$, occur.
In each case we compute the equality kernel in every degree and particle sector. 
The resulting decomposition gives the complete equality kernel for exclusion and, through the regular representation, for the interchange process.
\end{abstract}

\maketitle
\setcounter{tocdepth}{1}
\tableofcontents

\section{Introduction and main results}
\label{sec:introduction}

Let $G=(V,c)$ be a connected undirected weighted graph on $n$ vertices, with symmetric
conductances $c_{uv}=c_{vu}\ge0$, $c_{uu}=0$, for all $u,v\in V$.
We regard $\{u,v\}$ as an edge exactly when $c_{uv}>0$.  
Fix distinct terminals $x,y\in V$.
The \emph{effective conductance} between $x$ and $y$ is
\begin{equation}\label{eq:effective-conductance}
  \ceff(x,y)
  :=\inf\left\{
    \sum_{\{u,v\}}c_{uv}[\phi(u)-\phi(v)]^2 ~\bigg|~
    \phi:V\to\mathbb R,\ \phi(x)=0,\ \phi(y)=1
  \right\}.
\end{equation}
The energy functional in \eqref{eq:effective-conductance} has a unique minimizer, which we denote by $p$ and call the \emph{harmonic $x$--$y$ voltage}: it satisfies $p(x)=0$, $p(y)=1$, and the discrete Laplace's equation
\begin{equation}\label{eq:harmonic-voltage-equation}
  \sum_{v\in V\setminus \{u\}}c_{uv}[p(u)-p(v)]=0,
  \qquad u\in V\setminus\{x,y\}.
\end{equation}
The discrete maximum principle gives $0\le p\le1$ on $V$.
By the homogeneity of the variational characterization in
\eqref{eq:effective-conductance}, every $\phi:V\to\mathbb R$ satisfies
\begin{align}
\label{eq:Dirichlet-ineq}
  \ceff(x,y)[\phi(x)-\phi(y)]^2
  \le
  \sum_{\{u,v\}}c_{uv}[\phi(u)-\phi(v)]^2.
\end{align}

The \emph{moving particle lemma} (or inequality) is the interacting particle analog of the
inequality \eqref{eq:Dirichlet-ineq}: the scalar edge energy on the right-hand side is replaced by the Dirichlet form of the
particle system, while the terminal difference on the left-hand side is replaced by the energy cost of
exchanging the configuration at $x$ and $y$.

The moving particle lemma has a long history in the interacting particle systems literature; see, for example, Kipnis and
Landim~\cite[\S5.5]{KipnisLandim}.
In the classical argument, one chooses a path joining the two sites, writes the distant
exchange as a telescoping sequence of nearest-neighbor exchanges along that path, and
applies Cauchy--Schwarz.  The resulting constants are controlled by path lengths and
congestion.
Alternatively, one may replace this pathwise comparison by the electrical coefficient \eqref{eq:effective-conductance}.

Throughout this paper we work with the interchange process on $G$.  
Let $\mathfrak S_V$ denote the symmetric group on $V$; it acts by left multiplication in its regular representation on $\mathbb R[\mathfrak S_V]$.
For distinct $u,v\in V$, denote their transposition by $\tau_{uv}$.
We define the positive semidefinite interchange process generator as
\[
  \cL_G^{\mathrm{IP}}
  :=\sum_{\{u,v\}\subset V}c_{uv}(I-\tau_{uv}).
\]
In \cite{ChenMPL}, this author used the octopus inequality of Caputo, Liggett, and Richthammer
\cite{CLR}---the key ingredient towards resolving Aldous' spectral gap conjecture---to prove the following operator inequality.
\begin{proposition}[Moving particle lemma, \emph{cf.\@} Lemma 3.2 in \cite{ChenMPL}]
\label{prop:MPL}
On every finite connected weighted graph $G=(V,c)$, and distinct vertices $x,y\in V$,
\begin{equation}\label{eq:IP-MPL-intro}
  \cL_G^{\mathrm{IP}}
  \ge \ceff(x,y)\bigl(I-\tau_{xy}\bigr)
\end{equation}
as operators on the regular representation of $\mathfrak S_V$.
\end{proposition}
The proof of \cref{prop:MPL} proceeds via eliminating vertices successively while preserving the effective conductance between $x$ and $y$, and applying the octopus inequality \cite[Theorem 2.3]{CLR} at each step of the elimination.
The argument goes through in only a few lines, but hides the nontrivial nature of the octopus inequality, and does not utilize any algebraic aspect thereof. For a representation theory perspective of the octopus inequality, see \emph{e.g.\@} Cesi \cite{Cesi}.

Alon and Kozma \cite{AK20} later described
\cite{ChenMPL} as the first application of the octopus inequality for a purpose beyond the
proof of Aldous' spectral gap conjecture. 
Of course there have been many more developments since then, thanks to our better understanding of the algebraic structure of the operators involved.
The main focus of this work is to establish a full characterization of attaining equality in \eqref{eq:IP-MPL-intro}.

In order to proceed, let us discuss how to transfer \cref{prop:MPL} to the symmetric exclusion process on $G$.
For $0\le r\le n$, write $\Omega_r=\binomset Vr$ for the space of $r$-subsets (or the
$r$-particle sector in the exclusion process language), and equip $\mathbb R^{\Omega_r}$ with the normalized inner product
\begin{align}
  \langle f,g\rangle_r
  =\binom nr^{-1}\sum_{\xi\in\Omega_r}f(\xi)g(\xi),
  \qquad
  \|f\|_r^2=\langle f,f\rangle_r.
 \label{eq:r-inner-prod}
\end{align}
The positive semidefinite exclusion process generator is
\begin{equation}\label{eq:positive-generator}
  \cL_G^{(r)}
  =\sum_{\{u,v\}\subset V}c_{uv}(I-\tau_{uv}),
\end{equation}
where the transposition $\tau_{uv}$ acts on subsets by interchanging $u$ and $v$, and on
functions by pullback.  
Now apply \cref{prop:MPL} to the permutation representation on $r$-subsets to obtain
\begin{align}
\label{eq:MPL-exclusion}
\cL_G^{(r)}\ge \ceff(x,y)(I-\tau_{xy}).
\end{align}
It is standard to express \eqref{eq:MPL-exclusion} as a quadratic form inequality under the inner product \eqref{eq:r-inner-prod}; in particular at $r=1$ one recovers \eqref{eq:Dirichlet-ineq}.  
There is also a related product Bernoulli formulation of this quadratic form inequality that is
recorded as \cite[Theorem~1.1]{ChenMPL}.
For most of this paper, we focus on the equality analysis in each fixed $r$-particle sector.

This requires us to analyze how the operators act on the irreducible representations of $\mathfrak{S}_V$, or the so-called Specht modules.
For background and standard facts on Specht modules, Young diagrams, and the branching rule used below, we refer the reader to Sagan \cite{Sagan}.
For a given partition $\lambda=(\lambda_1, \cdots, \lambda_m)$ of $n$, write $S^\lambda$ for the Specht
module (over the field $\mathbb{R}$) of $\mathfrak S_V$ indexed by $\lambda$.  For $r\le \frac{n}{2}$,
Young's rule \cite[\S2.11]{Sagan} states that the $r$-subset permutation representation decomposes multiplicity-free as
\begin{align}
\label{eq:perm-Young}
  \mathbb R^{\Omega_r}
  \cong S^{(n)}\oplus S^{(n-1,1)}\oplus\cdots\oplus S^{(n-r,r)}.
\end{align}
This decomposition agrees with the Hoeffding decomposition to be used from \cref{sec:interface} onwards: the constant
sector is $S^{(n)}$, the linear sector is $S^{(n-1,1)}$, the quadratic sector is $S^{(n-2,2)}$, and, more generally, the
degree-$j$ Hoeffding layer is the copy of $S^{(n-j,j)}$.  
Thus symmetric exclusion realizes the interchange generator on the two-row irreducible representations in \eqref{eq:perm-Young}, with exactly the same conductance coefficients.

In \cite{ChenMPL} this author conjectured that the coefficient $\ceff(x,y)$ in \eqref{eq:IP-MPL-intro} should be sharp.
Recently Zhu \cite{Zhu26} proved the said conjecture on every $G$ by exhibiting a vector in the linear sector $S^{(n-1,1)}$ that realizes equality in \eqref{eq:IP-MPL-intro}; see also \cref{prop:degree-one-kernel} below.
This naturally leads to the next question studied in this work: which vectors outside the sectors $S^{(n)}$ and $S^{(n-1,1)}$ attain equality in \eqref{eq:IP-MPL-intro}?  We call any nonzero vector that attains equality in the relevant moving particle inequality an \emph{equality vector}.

A related but distinct problem concerning the second eigenspace of the interchange process generator
$\cL_G^{\mathrm{IP}}$ was studied by Belotserkovskiy and this author in \cite{BelotserkovskiyChen}; see \cref{rem:second-eigenspace-comparison} for the precise comparison.

\subsection{Summary of the main results}

The regular representation decomposes orthogonally as
\begin{equation}\label{eq:regular-specht-decomposition}
  \RR[\mathfrak S_V]
  \cong\bigoplus_{\lambda\vdash n}S^\lambda\otimes\RR^{d_\lambda},
  \qquad d_\lambda:=\dim S^\lambda,
\end{equation}
where each Specht module $S^\lambda$ occurs with multiplicity $d_\lambda$; see
\cite[\S1.10 and \S2.3]{Sagan}.  
Because the operators on both sides of \eqref{eq:IP-MPL-intro} are built
from permutations, the equality problem can be studied on each irreducible
summand.  
\Cref{thm:interchange-two-row-rigidity}
shows that  no nonzero equality vector can be supported on a partition with three or more rows.  
Consequently, all equality phenomena for the interchange process are already visible in
the two-row representations realized by symmetric exclusion.
We name this finding \emph{two-row rigidity}.

By carrying out the analysis in each $r$-particle sector of symmetric exclusion, we obtain a complete graph-theoretic classification, stated formally in \cref{thm:equality-beyond-degree-one}.  
Put $W:=V\setminus\{x,y\}$, let $G[W]$ be the induced graph on $W$ that retains the edges of positive conductance inside $W$,
and let $p$ be the harmonic $x$--$y$ voltage from \eqref{eq:harmonic-voltage-equation}.  
Then equality in the moving particle lemma outside the
constant sector $S^{(n)}$ and the linear sector $S^{(n-1,1)}$ occurs if and only if one of the following alternatives
holds:
\begin{enumerate}[label=\textup{(\alph*)},wide]
\item $p(u)=\frac12$ for every $u\in W$, and $G[W]$ has at least two connected
      components;
\item $p$ is constant on every connected component of $G[W]$, and the component
      values are exactly $a$ and $1-a$ for some $0\le a<\frac12$, with both values
      occurring.
\end{enumerate}
No other weighted graph and terminal pair supports equality beyond the constant and
linear sectors.  \Cref{thm:top-degree-kernel-classification,thm:particle-sector-kernel} determine
the dimension and structure of every surviving equality space, first in each top Hoeffding degree
and then in every particle sector.  Finally,
\cref{thm:complete-interchange-kernel} assembles these spaces into the complete equality
kernel for the interchange process.

The rest of the paper is organized as follows.
The proof begins with the aforementioned two-row rigidity in \cref{sec:interchange-rigidity},
and the equality criterion in symmetric exclusion in \cref{sec:equality-criterion}. 
In \cref{sec:interface} we derive interface identities from an equality vector and give conditions for reconstructing one from interface data. We also show that a nonzero equality vector in a top Hoeffding space $\cH_k$, with $3\le k\le\lfloor \frac{n}{2}\rfloor$, forces the nonterminal harmonic voltage to take at most two distinct values.  
We then classify the quadratic kernel in
\cref{sec:quadratic-classification}, and reduce higher-degree equality to it in
\cref{sec:degree-lowering}.  Finally, \cref{sec:sector-transfer,sec:higher-kernels}
transfer and reconstruct the kernels in all particle sectors, first for symmetric exclusion and then for the interchange process.

\section{Interchange equality and two-row rigidity}
\label{sec:interchange-rigidity}

Let $\lambda\vdash n$ be a partition, let $\ell(\lambda)$ denote the number of rows in its Young diagram, and let $\lambda_i$ be the number of boxes in row $i$.
Fix an $\mathfrak S_V$-invariant inner product on $S^\lambda$, and write
$\rho_\lambda$ for the corresponding orthogonal representation.
Set
\begin{align}
\label{eq:Q-L-K}
  Q_{uv}^{\lambda}:=\mathrm{Id}_{S^\lambda}-\rho_\lambda(\tau_{uv}),
  \qquad
  \cL_G^{\lambda}:=\sum_{\{u,v\}\subset V}c_{uv}Q_{uv}^{\lambda},
  \qquad
  \cK_{G,xy}^{\lambda}:=\cL_G^{\lambda}-\ceff(x,y)Q_{xy}^{\lambda},
\end{align}
where $x\neq y \in V$.
Since the regular representation contains every Specht module, the moving particle lemma for the interchange process, \cref{prop:MPL}, implies $\cK_{G,xy}^{\lambda}\ge0$.

\subsection{One-vertex reduction}

Fix a nonterminal vertex $z\in V\setminus\{x,y\}$, put $V_z:=V\setminus\{z\}$, and write
$\mathfrak S_{V_z}$ for the subgroup fixing $z$.  Also write $c_i:=c_{iz}$ and $C_z:=\sum_{i\in V_z}c_i>0$.
When the underlying graph matters, we write $\ceff^H(x,y)$ for the effective conductance between $x$ and $y$ in a weighted graph $H$; the definition still follows from \eqref{eq:effective-conductance}.  

Upon reducing the graph $G$ at the vertex $z$, we obtain the Schur-reduced graph $G/z$ on $V_z$, which is connected and has conductances
\begin{equation}\label{eq:Schur-weights}
  c^{(z)}_{ij}:=c_{ij}+\frac{c_i c_j}{C_z},
  \qquad i\ne j\in V_z.
\end{equation}
We claim that this reduction preserves the effective conductance between the terminals $x$ and $y$.
For a given voltage function $\phi:V_z\to \RR$, put
\[
  \bar\phi_z:=\frac1{C_z}\sum_{i\in V_z}c_i\phi(i).
\]
If $t$ denotes the voltage value at $z$, then the Dirichlet energy contributed by edges connected to $z$ equals
\[
  \sum_{i\in V_z}c_i(t-\phi(i))^2
  =C_z(t-\bar\phi_z)^2
   +\frac1{C_z}\sum_{\{i,j\}\subseteq V_z}c_i c_j[\phi(i)-\phi(j)]^2.
\]
From the right-hand side expression, we see that this energy attains the minimum at $t=\bar\phi_z$, with value exactly the second term.  Adding the edges already contained in $V_z$
therefore replaces $c_{ij}$ by the reduced conductance $c^{(z)}_{ij}$ in \eqref{eq:Schur-weights}.
Since $z$ is not a terminal, the constraints at $x$ and $y$ are unchanged, and a
second minimization over the values on $V_z$ gives
\begin{equation}\label{eq:Schur-ceff}
  \ceff^{G/z}(x,y)=\ceff^G(x,y).
\end{equation}

We now analyze two positive semidefinite operators that are relevant to the problem.  
The operator $\cK_{G,xy}^{\lambda}$ introduced in \eqref{eq:Q-L-K} is referred to as the \emph{moving-particle residual}.  
The operator created by eliminating $z$ and the edges attached to $z$ is the \emph{one-vertex octopus remainder}
\begin{equation}\label{eq:octopus-remainder}
  \Delta_z^\lambda
  :=\sum_{i\in V_z}c_iQ_{iz}^{\lambda}
    -\frac1{C_z}\sum_{\{i,j\}\subseteq V_z}c_i c_jQ_{ij}^{\lambda}.
\end{equation}
To see how the two operators fit together, write the residual on $G$ as
\begin{align}\label{eq:K_G-lambda}
  \cK_{G,xy}^{\lambda}
  =\sum_{\{i,j\}\subseteq V_z}c_{ij}Q_{ij}^{\lambda}
   +\sum_{i\in V_z}c_iQ_{iz}^{\lambda}
   -\ceff^G(x,y)Q_{xy}^{\lambda}.
\end{align}
On the restricted $\mathfrak S_{V_z}$-representation $S^\lambda\downarrow_{\mathfrak S_{V_z}}$, define the reduced residual by
\begin{align}\label{eq:K_G-lambda-V_z}
  \cK_{G/z,xy}^{\lambda\downarrow V_z}
  :=\sum_{\{i,j\}\subseteq V_z}
      \left(c_{ij}+\frac{c_ic_j}{C_z}\right)Q_{ij}^{\lambda}
   -\ceff^{G/z}(x,y)Q_{xy}^{\lambda}.
\end{align}
Subtracting \eqref{eq:K_G-lambda-V_z} from \eqref{eq:K_G-lambda} and using \eqref{eq:Schur-ceff} leaves
exactly \eqref{eq:octopus-remainder}.  Hence
\begin{equation}\label{eq:Schur-residual-decomposition}
  \cK_{G,xy}^{\lambda}
  =\Delta_z^\lambda+\cK_{G/z,xy}^{\lambda\downarrow V_z}.
\end{equation}

According to the branching rule, \emph{cf.\@} \cite[\S2.8]{Sagan}, the restriction of the Specht module $S^\lambda$ to the subgroup $\mathfrak{S}_{V_z}$ decomposes multiplicity-free as
\begin{equation}\label{eq:specht-branching}
  S^\lambda\!\downarrow_{\mathfrak S_{V_z}}
  \cong\bigoplus_{\mu=\lambda-\square}S^\mu,
\end{equation}
where the sum runs over the removable boxes of $\lambda$.

\begin{lemma}[Positivity inputs for one-vertex reduction]
\label{lem:one-vertex-positivity}
For every partition $\lambda\vdash n$,
\[
  \Delta_z^\lambda\ge0
  \quad\text{and}\quad
  \cK_{G/z,xy}^{\lambda\downarrow V_z}\ge0,
\]
where the first inequality acts on $S^\lambda$, and the second on the restricted
$\mathfrak S_{V_z}$-representation $S^\lambda\!\downarrow_{\mathfrak S_{V_z}}$.
\end{lemma}

\begin{proof}
The first inequality is the octopus inequality of Caputo, Liggett, and Richthammer
\cite[Theorem 2.3]{CLR} restricted to the representation $S^\lambda$.  
For the second inequality, apply \cref{prop:MPL} to the connected
Schur-reduced graph $G/z$.  By the branching decomposition
\eqref{eq:specht-branching}, the restricted representation
$S^\lambda\!\downarrow_{\mathfrak S_{V_z}}$ is the direct sum of the irreducible
summands $S^\mu$ obtained by removing one box from $\lambda$.  The resulting
operator inequality is nonnegative on each such summand, and hence on their direct
sum.
\end{proof}

Let us make a simple but crucial observation.
If $\ell(\lambda)\ge 3$, then every branch in \eqref{eq:specht-branching} still has at least 3 rows except in exactly one scenario: $\lambda=(a,b,1)$ with $a\ge b\ge 1$, and the box in row 3 is removed, producing $\mu=(a,b)$.
We call this the \emph{exceptional branch}.

To prepare for the next lemma, recall that if $P$ is the orthogonal projection onto a
branching summand, then $P\Delta_z^\lambda P$ is the compression of
$\Delta_z^\lambda$ to that summand; equivalently, it is the operator governing the
quadratic form of $\Delta_z^\lambda$ on vectors in the range of $P$.  We also use the
standard Jucys--Murphy content rule for branching (see Murphy~\cite{Murphy81}).
In a Young diagram, the \emph{content} of a box in row $r$ and column $c$ is
$\operatorname{ct}(r,c)=c-r$.  If a branch
$S^\mu\subset S^\lambda\!\downarrow_{\mathfrak S_{V_z}}$ is obtained by removing the box
$\lambda\setminus\mu$, and $P_\mu$ is the orthogonal projection onto that branch, then the Jucys--Murphy element $\sum_{i\in V_z}\rho_\lambda(\tau_{iz})$ acts as a scalar on $P_\mu$ according to
\begin{equation}\label{eq:JM-branch-content}
  \left(\sum_{i\in V_z}\rho_\lambda(\tau_{iz})\right)P_\mu
  =\operatorname{ct}(\lambda\setminus\mu)P_\mu.
\end{equation}
The sum in \eqref{eq:JM-branch-content} is the
difference between the transposition class sums for $\mathfrak S_V$ and
$\mathfrak S_{V_z}$, and these act on $S^\lambda$ and $S^\mu$ by the respective sums of box
contents.  

For the said exceptional branch, fix $\lambda=(a,b,1)$ with $a\ge b\ge1$, let
$B_z\cong S^{(a,b)}$ be the branch obtained by removing the third-row box, and let $P$
be the orthogonal projection onto $B_z$.  In the next lemma,
we write $Q_{uv}:=Q_{uv}^{(a,b,1)}$; restrictions of these operators to $B_z$
are indicated explicitly.

\begin{lemma}[Exceptional branch compression]
\label{lem:exceptional-branch-compression}
With the notation above, the following hold:
\begin{enumerate}[label=\textup{(\roman*)},ref=(\roman*)]
\item\label{it:exceptional-compressed-transposition}
For every $i\in V_z$,
\[
  \left(PQ_{iz}P\right)\big|_{B_z}
  =\mathrm{Id}_{B_z}
   +\frac1{b(a+1)}
      \sum_{j\in V_z\setminus\{i\}}Q_{ij}\big|_{B_z}.
\]

\item\label{it:exceptional-positivity}
The octopus remainder satisfies
\[
  P\Delta_z^{(a,b,1)}P
  \ge \frac{\sum_{i\in V_z}c_i^2}{C_z}\,P,
\]
and hence is strictly positive on $B_z$.
\end{enumerate}
\end{lemma}

\begin{proof}
Write $\rho:=\rho_{(a,b,1)}$ and regard all compressed operators as acting on $B_z$.  The removed box is $(3,1)$, whose content is $-2$, so \eqref{eq:JM-branch-content} gives
\[
  \left(\sum_{j\in V_z}\rho(\tau_{jz})\right)\bigg|_{B_z}
  =-2\,\mathrm{Id}_{B_z}.
\]
Fix $i\in V_z$.  The group algebra identity
\[
  \tau_{iz}\sum_{j\in V_z}\tau_{jz}
  =1+\left(\sum_{j\in V_z\setminus\{i\}}\tau_{ij}\right)\tau_{iz}
\]
may be compressed to $B_z$ as follows.  The projection $P$ commutes with
every $\rho(\tau_{ij})$, $i,j\in V_z$, since $B_z$ is an invariant summand of the
orthogonal subgroup representation.  Also, \eqref{eq:JM-branch-content} gives
$(\sum_{j\in V_z}\rho(\tau_{jz}))P=-2P$.  Multiplying by $P$ on both sides gives
\[
  \left(
    2\,\mathrm{Id}_{B_z}
    +\sum_{j\in V_z\setminus\{i\}}\rho(\tau_{ij})\big|_{B_z}
  \right)
  \left(P\rho(\tau_{iz})P\right)\big|_{B_z}
  =-\mathrm{Id}_{B_z}.
\]
After relabeling $i$ as the last letter of $V_z$, the sum over $j$ is the Jucys--Murphy element for restricting $S^{(a,b)}$ to $\mathfrak S_{V_z\setminus\{i\}}$.  Its possible eigenvalues are the contents of the removable boxes of $(a,b)$, hence lie in $\{a-1,b-2\}$.  For either such eigenvalue $t$,
$
  (2+t)^{-1}=\frac{a+b-1-t}{b(a+1)}$.
Spectral calculus therefore gives
\[
  \left(
    2\,\mathrm{Id}_{B_z}
    +\sum_{j\in V_z\setminus\{i\}}\rho(\tau_{ij})\big|_{B_z}
  \right)^{-1}
  =\frac1{b(a+1)}
   \sum_{j\in V_z\setminus\{i\}}Q_{ij}\big|_{B_z},
\]
which proves item~\ref{it:exceptional-compressed-transposition}.

For item~\ref{it:exceptional-positivity}, substitute item~\ref{it:exceptional-compressed-transposition} into \eqref{eq:octopus-remainder}.  This gives
\[
  \left(P\Delta_zP\right)\big|_{B_z}
  =C_z\,\mathrm{Id}_{B_z}
   +\frac1{b(a+1)}
      \sum_{i\in V_z}c_i
      \sum_{j\in V_z\setminus\{i\}}Q_{ij}\big|_{B_z}
  -\frac1{C_z}\sum_{\{i,j\}\subseteq V_z}c_i c_jQ_{ij}\big|_{B_z}.
\]
The double sum is positive semidefinite, and each transposition is a self-adjoint involution, so $0\le Q_{ij}|_{B_z}\le2\,\mathrm{Id}_{B_z}$.  Hence
\[
  \left(P\Delta_zP\right)\big|_{B_z}
  \ge
  \left(C_z-\frac2{C_z}\sum_{\{i,j\}\subseteq V_z}c_i c_j\right)\mathrm{Id}_{B_z}
  =\frac{\sum_{i\in V_z}c_i^2}{C_z}\,\mathrm{Id}_{B_z},
\]
using $C_z^2=\sum_i c_i^2+2\sum_{\{i,j\}\subseteq V_z}c_ic_j$.  This proves item~\ref{it:exceptional-positivity}.
\end{proof}

\subsection{No equality in representations with three or more rows}

The next result says that equality in the moving particle lemma can occur only in irreducibles with at most two rows.

\begin{theorem}[Two-row rigidity for interchange equality]
\label{thm:interchange-two-row-rigidity}
Let $G=(V,c)$ be a finite connected weighted graph, let $x\ne y\in V$ be terminals, and let
$\lambda\vdash n=|V|$.  If $\ell(\lambda)\ge3$, then
$\ker\cK_{G,xy}^{\lambda}=\{0\}$.
\end{theorem}

\begin{proof}
We argue by induction on $n$.  The assertion is vacuous for $n\le2$.  Assume it is
known on every graph with $n-1$ vertices, and let $\ell(\lambda)\ge3$.  Choose
$z\notin\{x,y\}$ and suppose that $0\ne f\in\ker\cK_{G,xy}^{\lambda}$.  By
\eqref{eq:Schur-residual-decomposition} and \cref{lem:one-vertex-positivity},
\begin{equation}\label{eq:equality-splits-under-Schur}
  \Delta_z^\lambda f=0
  \quad\text{and}\quad
  \cK_{G/z,xy}^{\lambda\downarrow V_z}f=0.
\end{equation}
Decompose $f$ according to the branching rule \eqref{eq:specht-branching}.  The second
operator in \eqref{eq:equality-splits-under-Schur} belongs to the group algebra of
$\mathfrak S_{V_z}$, so it preserves every branch; in particular on the branch $S^\mu$ it is $\cK_{G/z,xy}^{\mu}$.  
By the induction hypothesis, its kernel is trivial whenever $\ell(\mu)\ge3$.

If no removable box of $\lambda$ produces a two-row partition, every branch component
of $f$ is therefore zero, a contradiction.  The only remaining possibility is
$\lambda=(a,b,1)$, and then the only two-row branch is the branch $B_z\cong S^{(a,b)}$
of \cref{lem:exceptional-branch-compression}.  Hence $f\in B_z$, so $Pf=f$.  The first
identity in \eqref{eq:equality-splits-under-Schur} and
\cref{lem:exceptional-branch-compression}\ref{it:exceptional-positivity} give
\[
  0=\langle f,\Delta_z^\lambda f\rangle
   =\langle f,P\Delta_z^\lambda P f\rangle
   \ge \frac{\sum_{i\in V_z} c_i^2}{C_z}\|f\|^2>0,
\]
a contradiction.  This completes the induction.
\end{proof}

\begin{remark}
\label{rem:second-eigenspace-comparison}
The proof of \cref{thm:interchange-two-row-rigidity} shares some
representation theoretic ingredients with the proof strategy in \cite{BelotserkovskiyChen}.
Both works track the equality of operators
through a one-vertex reduction, and use the branching rule to control the irreducible
components that survive restriction to one fewer vertex.  
But the objects of study differ.  
In \cite{BelotserkovskiyChen} the goal is to identify the 2nd eigenspace of $\mathcal{L}^{\rm IP}_G$.  
In this work, we fix the terminals $x,y$, and study the kernel of the terminal-dependent operator
$\cK_{G,xy}^{\lambda}$.
\end{remark}

\section{Equality criterion for symmetric exclusion}
\label{sec:equality-criterion}

Having ruled out partitions with three or more rows, we now analyze the remaining
two-row representations through symmetric exclusion in an arbitrary $r$-particle sector.

Retain the notation of \cref{sec:introduction}.  Since each $\tau_{uv}$ is an
orthogonal involution, the quadratic form of \eqref{eq:positive-generator} is
\begin{equation}\label{eq:exclusion-energy}
  \langle f,\cL_G^{(r)}f\rangle_r
  =\sum_{\{u,v\}}c_{uv}\langle f,(I-\tau_{uv})f\rangle_r, \qquad f\in \mathbb{R}^{\Omega_r}.
\end{equation}
Recall from \eqref{eq:harmonic-voltage-equation} that the minimizing voltage $p$
is harmonic at every nonterminal vertex.  
In what follows we will often use the concatenation notation:
if $u\notin T$, then $uT$ means $T\sqcup\{u\}$, and similarly for $uvT$.

\subsection{Transposition metric and the equality criterion}

Throughout this subsection $1\le r\le n-1$; the sectors $r=0,n$ consist only of
constants.  The key geometric input is that transposition differences of a fixed
observable $f\in \mathbb R^{\Omega_r}$ define a pseudometric on the vertex set $V$.

\begin{definition}[Transposition distance]\label{def:transposition-distance}
For $f\in\RR^{\Omega_r}$ and distinct $u,v\in V$, define
\begin{equation}\label{eq:transposition-distance}
  \delta_f(u,v)^2:=\langle f,(I-\tau_{uv})f\rangle_r
  =\frac12\|f-\tau_{uv}f\|_r^2.
\end{equation}
Set $\delta_f(u,u):=0$ for every $u\in V$.
For distinct $u,v\in V$, the same quantity has the coordinate representation
\begin{equation}\label{eq:delta-coordinate}
  \delta_f(u,v)^2
  =\binom nr^{-1}
    \sum_{T\in\binomset{V\setminus\{u,v\}}{r-1}}
    \bigl[f(uT)-f(vT)\bigr]^2.
\end{equation}
We write
\[
  \mathbf d_f^{uv}
  :=\bigl(f(uT)-f(vT)\bigr)_{T\in\binomset{V\setminus\{u,v\}}{r-1}}
\]
for this coordinate difference vector.  For distinct $u,v,w\in V$, define two further coordinate vectors indexed by the same sets $T$ by
\begin{align}
\label{eq:abuvfw}
\bigl(a_{f,w}^{uv}(T),b_{f,w}^{uv}(T)\bigr):=
\begin{cases}
\bigl(f(uT)-f(wT),\ f(wT)-f(vT)\bigr),&w\notin T,\\
\bigl(f(uvT')-f(vwT'),\ f(uwT')-f(uvT')\bigr),&T=T'\sqcup\{w\},
\end{cases}
\end{align}
and write $\mathbf a_{f,w}^{uv}:=(a_{f,w}^{uv}(T))_{T\in\binomset{V\setminus\{u,v\}}{r-1}}$ and
$\mathbf b_{f,w}^{uv}:=(b_{f,w}^{uv}(T))_{T\in\binomset{V\setminus\{u,v\}}{r-1}}$.
\end{definition}

\begin{lemma}[Transposition pseudometric and equality in the triangle inequality]
\label{lem:transposition-metric}
For every $f\in\RR^{\Omega_r}$, the following hold for distinct $u,v,w\in V$:
\begin{enumerate}[label=\textup{(\roman*)},ref=(\roman*)]
\item\label{it:transposition-splitting}
The coordinate vectors split as $\mathbf d_f^{uv}=\mathbf a_{f,w}^{uv}+\mathbf b_{f,w}^{uv}$,
and their Euclidean norms are
\[
  \|\mathbf d_f^{uv}\|_2=\binom nr^{\frac12}\delta_f(u,v),\qquad
  \|\mathbf a_{f,w}^{uv}\|_2=\binom nr^{\frac12}\delta_f(u,w),\qquad
  \|\mathbf b_{f,w}^{uv}\|_2=\binom nr^{\frac12}\delta_f(w,v).
\]

\item\label{it:transposition-triangle}
The triangle inequality $\delta_f(u,v)\le\delta_f(u,w)+\delta_f(w,v)$ holds.
Hence $\delta_f$ is a pseudometric on $V$.

\item\label{it:transposition-equality}
If equality holds in item~\ref{it:transposition-triangle} and
$\delta_f(u,v)>0$, then
\[
  \mathbf a_{f,w}^{uv}
  =\frac{\delta_f(u,w)}{\delta_f(u,v)}\,\mathbf d_f^{uv},
  \qquad
  \mathbf b_{f,w}^{uv}
  =\frac{\delta_f(w,v)}{\delta_f(u,v)}\,\mathbf d_f^{uv}.
\]
\end{enumerate}
\end{lemma}

\begin{proof}
Item~\ref{it:transposition-splitting}: The vector identity holds coordinatewise.
When $w\notin T$, $a_{f,w}^{uv}(T)$ is a $(u,w)$-difference and $b_{f,w}^{uv}(T)$ is a
$(w,v)$-difference.  When $T=T'\sqcup\{w\}$, the same is true after reindexing by
$vT'$ and $uT'$, respectively.  The norm identities then follow from
\eqref{eq:delta-coordinate}.

Item~\ref{it:transposition-triangle}: Symmetry follows from
\eqref{eq:transposition-distance}, and $\delta_f(u,u)=0$ by definition.
Minkowski's inequality applied to
item~\ref{it:transposition-splitting} gives the triangle inequality.

Item~\ref{it:transposition-equality}: Equality in
item~\ref{it:transposition-triangle} is equality in the Euclidean triangle
inequality from item~\ref{it:transposition-splitting}.  Hence
$\mathbf a_{f,w}^{uv}$ and $\mathbf b_{f,w}^{uv}$ are nonnegative scalar multiples of
$\mathbf d_f^{uv}$.
Comparing their norms in item~\ref{it:transposition-splitting}
gives the stated factors.
\end{proof}

\begin{proposition}[Moving particle lemma and equality criterion]
\label{prop:MPL-equality}
For every $1\le r\le n-1$, every $f\in\RR^{\Omega_r}$, and every $x\neq y$,
\begin{equation}\label{eq:MPL-direct}
  \ceff(x,y)\,\delta_f(x,y)^2
  \le \langle f,\cL_G^{(r)}f\rangle_r.
\end{equation}
If $\delta_f(x,y)>0$, equality holds in \eqref{eq:MPL-direct} if and only if the function
$\phi_f(v):=\delta_f(x,v)/\delta_f(x,y)$ is the harmonic $x$--$y$ voltage and every
edge $\{u,v\}$ of positive conductance satisfies
\begin{equation}\label{eq:edge-equality}
  \delta_f(u,v)
  =\delta_f(x,y)|\phi_f(u)-\phi_f(v)|.
\end{equation}
If $\delta_f(x,y)=0$, equality holds if and only if $f$ is constant on $\Omega_r$.
\end{proposition}

\begin{proof}
Assume first that $\delta_f(x,y)>0$.  The reverse triangle inequality for the
pseudometric gives
$|\phi_f(u)-\phi_f(v)|\le \delta_f(u,v)/\delta_f(x,y)$.
Since $\phi_f(x)=0$ and $\phi_f(y)=1$, the scalar Dirichlet principle
\eqref{eq:effective-conductance} yields
\begin{equation}\label{eq:MPL-equality-chain}
  \ceff(x,y)
  \le \sum_{\{u,v\}}c_{uv}[\phi_f(u)-\phi_f(v)]^2
  \le \delta_f(x,y)^{-2}
       \sum_{\{u,v\}}c_{uv}\delta_f(u,v)^2.
\end{equation}
By \eqref{eq:exclusion-energy}, the sum over $\{u,v\}$ in the last display is
$\langle f,\cL_G^{(r)}f\rangle_r$, which then proves \eqref{eq:MPL-direct}.

Equality in \eqref{eq:MPL-direct} holds exactly when both inequalities in
\eqref{eq:MPL-equality-chain} are equalities.  To attain equality in the first inequality, $\phi_f$ must be the harmonic
voltage per \eqref{eq:effective-conductance}--\eqref{eq:harmonic-voltage-equation}.  
In the second inequality, every edge contribution is nonnegative, so equality
is equivalent to \eqref{eq:edge-equality} on every edge of positive conductance.

Now assume $\delta_f(x,y)=0$. 
Equality in \eqref{eq:MPL-direct} means
$\langle f,\cL_G^{(r)}f\rangle_r=0$.  Hence $\tau_{uv}f=f$ on every edge of positive
conductance.  The edge transpositions of a connected graph generate the full symmetric
group on $V$, so $f$ is constant on the $r$-particle sector.  The converse is immediate.
\end{proof}

From now on it is convenient to encode equality in the moving particle lemma by the residual operator
\begin{equation}\label{eq:exclusion-residual}
  \cK_{xy}^{(r)}
  :=\cL_G^{(r)}-\ceff(x,y)(I-\tau_{xy}).
\end{equation}

\begin{corollary}[Residual kernel and equality vectors]\label{cor:operator-equality}
For every $r$, $\cK_{xy}^{(r)}\ge0$.  Also, for $f\in\RR^{\Omega_r}$ we have
\[
  \cK_{xy}^{(r)}f=0
  \quad\Longleftrightarrow\quad
  \langle f,\cK_{xy}^{(r)}f\rangle_r=0
  \quad\Longleftrightarrow\quad
  \ceff(x,y)\delta_f(x,y)^2
  =\langle f,\cL_G^{(r)}f\rangle_r.
\]
If any of these conditions holds for a nonconstant $f$, then $\delta_f(x,y)>0$.
\end{corollary}

\begin{proof}
By \cref{prop:MPL}, applied in the $r$-subset permutation representation,
$\cL_G^{(r)}\ge \ceff(x,y)(I-\tau_{xy})$, so
\eqref{eq:exclusion-residual} gives $\cK_{xy}^{(r)}\ge0$.  For a positive
semidefinite self-adjoint operator, vanishing of the quadratic form is equivalent to
kernel membership.  Expanding the quadratic form of \eqref{eq:exclusion-residual}
gives the third condition.  Finally, \cref{prop:MPL-equality} shows that equality with
$\delta_f(x,y)=0$ forces $f$ to be constant.
\end{proof}

\subsection{Terminal geometry of equality vectors}

Fix a nonconstant equality vector $f$.  Write $u\sim_f v$ when
$\delta_f(u,v)=0$.  Since $\delta_f$ is a pseudometric, $\sim_f$ is an equivalence
relation on $V$; its equivalence classes are the \emph{zero-distance classes} of $f$.
Contracting these classes gives a quotient graph whose vertices are the zero-distance
classes and whose positive edges are induced by edges of positive conductance in $G$.

The next lemma turns the equality criterion into terminal geodesicity, and then invokes \cref{lem:transposition-metric}\ref{it:transposition-equality} to obtain the coordinate proportionality used in \cref{sec:interface}.

\begin{lemma}[Every vertex lies on a terminal geodesic]
\label{lem:every-vertex-geodesic}
Suppose $f$ is a nonconstant equality vector. 
Then for every $z\in V$,
\begin{equation}\label{eq:terminal-geodesic}
  \delta_f(x,y)=\delta_f(x,z)+\delta_f(z,y).
\end{equation}
Let $p_v:=p(v)=\phi_f(v)$, where $\phi_f$ is the harmonic $x$--$y$ voltage from \cref{prop:MPL-equality}.
Then for every $z\in V\setminus\{x,y\}$ we have
\begin{equation}\label{eq:terminal-coordinate-proportionality}
  \mathbf a_{f,z}^{xy}=p_z\mathbf d_f^{xy}
  \quad\text{and}\quad
  \mathbf b_{f,z}^{xy}=(1-p_z)\mathbf d_f^{xy},
\end{equation}
where $\mathbf a_{f,z}^{xy}$ and $\mathbf b_{f,z}^{xy}$ are the coordinate vectors from \eqref{eq:abuvfw} in  \cref{def:transposition-distance}.
\end{lemma}

\begin{proof}
Throughout the proof we use the shorthand $L_f := \delta_f(x,y)$.
By \cref{cor:operator-equality}, $L_f>0$.
For $z=x$ or $z=y$, \eqref{eq:terminal-geodesic} is immediate, so fix
$z\in V\setminus\{x,y\}$.  Pass to the quotient by the zero-distance relation $\sim_f$; the pseudometric $\delta_f$
descends to a metric on the quotient.  
Since $p_v=\delta_f(x,v)/L_f$, the reverse triangle inequality gives
$|p_u-p_v|\le\delta_f(u,v)/L_f$.  Thus $p$ is constant on each
zero-distance class.  Moreover, if distinct quotient classes $U,V$ are joined by an edge
induced from an edge of positive conductance in $G$, then
$ \delta_f(U,V)=L_f\,|p(U)-p(V)|$.
Thus adjacent quotient classes have different voltages, and along any quotient path
on which the voltage is monotone, the $\delta_f$-length is $L_f$ times the total
voltage change.

If the zero-distance class of $z$ contains $x$, then
$\delta_f(x,z)=0$, and the triangle inequality gives
$\delta_f(z,y)=\delta_f(x,y)$; the case in which it contains $y$ is
symmetric.  Thus \eqref{eq:terminal-geodesic} already holds in these cases.
We may then assume that the class $Q$ of $z$ contains neither terminal.
Choose a vertex of $Q$ incident to an edge
leaving $Q$, which exists because the quotient graph is connected.  If $p(Q)=0$, every
term in the harmonicity equation \eqref{eq:harmonic-voltage-equation} at that vertex
coming from outside $Q$ is nonpositive, and at least one is strictly negative, a
contradiction. The case $p(Q)=1$ is analogous.  Thus $0<p(Q)<1$.  The same harmonicity
argument shows that $Q$ has both a lower-voltage and a higher-voltage neighboring
class.  Because the quotient is finite, and voltage strictly decreases or increases at every
such step, a path obtained by following lower- or higher-voltage neighbors must
terminate.  No class without a terminal can be an endpoint: it has neighbors of
both voltage signs.  The maximum principle gives $0\le p\le1$, so the only
possible endpoints are respectively the class containing $x$ (voltage $0$) and
the class containing $y$ (voltage $1$).

Let $Q_0,\ldots,Q_m$ be such a decreasing-voltage path from the $z$-class to the
$x$-class.  Then
\[
  \sum_{i=0}^{m-1}\delta_f(Q_i,Q_{i+1})
  =L_f\sum_{i=0}^{m-1}\bigl(p(Q_i)-p(Q_{i+1})\bigr)
  =p_zL_f.
\]
Similarly, an increasing-voltage path from the $z$-class to the $y$-class has
$\delta_f$-length $(1-p_z)L_f$.  Hence
\[
  \delta_f(x,z)\le p_zL_f,
  \qquad
  \delta_f(z,y)\le (1-p_z)L_f.
\]
Together with the triangle inequality,
\[
  L_f=\delta_f(x,y)
  \le \delta_f(x,z)+\delta_f(z,y)
  \le p_zL_f+(1-p_z)L_f=L_f.
\]
Equality therefore holds throughout.  This proves \eqref{eq:terminal-geodesic}, and
\cref{lem:transposition-metric}\ref{it:transposition-equality} applied to the triple $(u,w,v)=(x,z,y)$ gives
\eqref{eq:terminal-coordinate-proportionality}.
\end{proof}

\section{Interface identities and voltage rigidity}
\label{sec:interface}

We introduce the \emph{down map} on an arbitrary ground set.  
Given a finite $E\subseteq V$ and $1\le r\le |E|$, define
\[
  D_r^E:\RR^{\binomset Er}\longrightarrow\RR^{\binomset E{r-1}},
  \qquad (D_r^E g)(T):=\sum_{z\in E\setminus T}g(T\cup\{z\}).
\]
When $E=V$, write $D_r:=D_r^V$ for short.  
Set
\[
  \cH_0(E):=\RR,\qquad
  \cH_r(E):=\ker D_r^E,
  \quad 1\le r\le\left\lfloor \frac{|E|}{2}\right\rfloor,
\]
and abbreviate $\cH_r:=\cH_r(V)$.  Thus $g\in\cH_r(E)$ if and only if $D_r^E g=0$.  The down map is equivariant under vertex permutations
of $E$, so $\cH_r(E)$ is preserved by the full symmetric group
$\mathfrak S_E$.  In particular, the space $\cH_r$ is
isomorphic to $S^{(n-r,r)}$ as an $\mathfrak S_V$-representation, as in
the Hoeffding--Specht correspondence of \cref{sec:introduction}; for
background on this realization see \cite{Filmus}.

\subsection{Interface identities and algebraic reconstruction}

This subsection provides the algebraic basis for our main theorems.
First, we derive in \cref{lem:interface-identities} identities for the terminal difference $h$ and the
both-terminal trace $b$ of a nonzero equality vector in $\cH_k$.
Then, in the reverse direction, we prove in \cref{lem:interface-reconstruction} that using compatible data $(p,h,b)$ one may reconstruct a function in $\cH_k$.
This reconstructed function becomes an equality vector under the harmonic voltage and the remaining positive-edge conditions, as stated in \cref{cor:reconstruction-equality}.

The six identities in \cref{lem:interface-identities} support several results to follow, such as the voltage rigidity result of \cref{sec:voltage-rigidity}, the degree lowering results in \cref{sec:degree-lowering}, and the classification of nonzero equality vectors in \cref{sec:higher-kernels}.

\begin{assumption}[Higher-degree equality hypothesis]\label{ass:top-degree-equality}
For the necessary interface identities below, the voltage-rigidity argument,
and the degree-lowering argument in \cref{sec:degree-lowering}, assume
\[
  3\le k\le \left\lfloor\frac n2\right\rfloor,
  \qquad
  0\neq f\in\cH_k,
  \qquad
  \cK_{xy}^{(k)}f=0.
\]
The function $f$ is nonconstant, since no nonzero constant belongs to $\cH_k$.
By \cref{prop:MPL-equality,cor:operator-equality},
$\delta_f(x,y)>0$ and $\phi_f=p$, where $p$ is the unique harmonic
$x$--$y$ voltage of $G$ (with $p(x)=0$ and $p(y)=1$).
Set $W:=V\setminus\{x,y\}$, $N:=|W|=n-2$, and $p_u:=p(u)$ for $u\in W$.
\end{assumption}

For any $A\subseteq W$, write $S_A:=\sum_{w\in A}p_w$; in particular,
$S_W$ is the total nonterminal voltage sum.
For $\xi\in\binomset W{k-1}$ and $T\in\binomset W{k-2}$, define
\begin{equation}\label{eq:interface-data}
  h_\xi:=f(x\xi)-f(y\xi),
  \qquad
  b_T:=f(xyT).
\end{equation}

\begin{lemma}[Interpolation and interface identities]\label{lem:interface-identities}
Under \cref{ass:top-degree-equality}, the interface functions $h$ and $b$ satisfy the following identities:
\begin{enumerate}[label=\textup{(\roman*)},ref=(\roman*)]
\item\label{it:interface-down} (Down-map identities)
$D_{k-1}^W h=0$ and $D_{k-2}^W b=0$.

\item\label{it:interface-terminal} (Terminal interpolation formula)
If $T\in\binomset W{k-2}$ and $z\in W\setminus T$, then
$f(xzT)=b_T+(1-p_z)h_{zT}$ and $f(yzT)=b_T-p_zh_{zT}$.

\item\label{it:interface-no-terminal} (Non-terminal interpolation formula)
If $\xi\in\binomset W{k-1}$ and $z\in W\setminus\xi$, then
$f(z\xi)=(1-p_z)f(x\xi)+p_zf(y\xi)$.

\item\label{it:interface-symmetry} (Symmetry)
If $R\in\binomset W{k-3}$ and $u,v\in W\setminus R$ are distinct, then
$b_{vR}-b_{uR}=(p_u-p_v)h_{uvR}$.

\item\label{it:interface-consistency} (Consistency)
If $T\in\binomset W{k-2}$ and $u,v\in W\setminus T$ are distinct, then
$(1-p_u-p_v)(h_{uT}-h_{vT})=0$.

\item\label{it:interface-master} (Master identity)
Writing $L:=N-k+2$, for every
$T\in\binomset W{k-2}$ and $u\in W\setminus T$ one has
\[
  (L+1)b_T+
  \left(L-S_W-Lp_u+S_T\right)h_{uT}=0.
\]
\end{enumerate}
\end{lemma}

\begin{proof}
Item~\ref{it:interface-down}: Let $R\in\binomset W{k-2}$.  Applying $D_kf=0$
to $xR$ gives $b_R+\sum_{z\in W\setminus R}f(xzR)=0$; applying it to $yR$
gives the same identity with $x$ replaced by $y$.  Subtraction proves $D_{k-1}^W h=0$.
Applying $D_kf=0$ to $xyQ$, with $Q\in\binomset W{k-3}$, gives $D_{k-2}^W b=0$.

Item~\ref{it:interface-terminal}: Apply the proportionality condition
\eqref{eq:terminal-coordinate-proportionality} of \cref{lem:every-vertex-geodesic} at the coordinate indexed by $zT$; note that
the $\mathbf a_{f,z}^{xy}$-coordinate is $b_T-f(yzT)$, the
$\mathbf b_{f,z}^{xy}$-coordinate is $f(xzT)-b_T$, and the
$\mathbf d_f^{xy}$-coordinate is $h_{zT}$.  

Item \ref{it:interface-no-terminal}: Apply the proportionality condition \eqref{eq:terminal-coordinate-proportionality} once more, now to $\xi$.

Item~\ref{it:interface-symmetry}: Fix $R\in\binomset W{k-3}$ and distinct
$u,v\in W\setminus R$.  Apply item~\ref{it:interface-terminal} to $f(xuvR)$,
first choosing $z=u$ and $T=vR$, and then choosing $z=v$ and $T=uR$.
This gives
\[
f(xuvR)
=b_{vR}+(1-p_u)h_{uvR}
=b_{uR}+(1-p_v)h_{uvR}.
\]
Manipulating the last equality yields $b_{vR}-b_{uR}=(p_u-p_v)h_{uvR}$, as required.

Item~\ref{it:interface-consistency}: Fix $T\in\binomset W{k-2}$ and distinct
$u,v\in W\setminus T$.  Apply item~\ref{it:interface-no-terminal} to $f(uvT)$,
taking $z=u$ and $\xi=vT$.  Substituting the two formulas from
item~\ref{it:interface-terminal}, we obtain
\[
\begin{aligned}
f(uvT)
&=(1-p_u)f(xvT)+p_uf(yvT)\\
&=(1-p_u)\bigl[b_T+(1-p_v)h_{vT}\bigr]
  +p_u\bigl[b_T-p_vh_{vT}\bigr]\\
&=b_T+(1-p_u-p_v)h_{vT}.
\end{aligned}
\]
Interchanging $u$ and $v$ in the same calculation gives
$
f(uvT)=b_T+(1-p_u-p_v)h_{uT}$.
Subtracting these two expressions gives
$(1-p_u-p_v)(h_{uT}-h_{vT})=0$.

Item~\ref{it:interface-master}: Fix $T\in\binomset W{k-2}$ and
$u\in W\setminus T$.  The parameter $L=N-k+2 \ge 2$ is convenient because $uT$
has $k-1$ elements, so $|W\setminus(uT)|=N-k+1=L-1 \ge 1$.  Thus $D_kf=0$ at
$uT$ consists of two terminal terms and $L-1$ nonterminal terms:
\begin{align}
\label{eq:fuT}
  0=f(xuT)+f(yuT)+\sum_{z\in W\setminus(uT)}f(zuT).
\end{align}
By item~\ref{it:interface-terminal}, the first two terms sum to
$2b_T+(1-2p_u)h_{uT}$.  For each $z\in W\setminus(uT)$,
item~\ref{it:interface-no-terminal} followed by item~\ref{it:interface-terminal} gives
\[
\begin{aligned}
  f(zuT)
  &=(1-p_z)f(xuT)+p_zf(yuT)\\
  &=(1-p_z)\bigl[b_T+(1-p_u)h_{uT}\bigr]
    +p_z\bigl[b_T-p_uh_{uT}\bigr]\\
  &=b_T+(1-p_u-p_z)h_{uT}.
\end{aligned}
\]
By the subset-sum convention, $S_W$ is the total nonterminal voltage sum.
The disjoint decomposition $W=\{u\}\sqcup T\sqcup\bigl(W\setminus(uT)\bigr)$
therefore gives $S_{W\setminus(uT)}=S_W-p_u-S_T$.
With these inputs \eqref{eq:fuT} rewrites as
\begin{align*}
0&=(L+1)b_T+
 \left(1-2p_u+(L-1)(1-p_u)
  -S_{W\setminus(uT)}\right)h_{uT}\\
 &=(L+1)b_T+ \left(L-(L+1)p_u-S_{W\setminus(uT)}\right)h_{uT}\\
 &=(L+1)b_T+
 \left(L-S_W-Lp_u+S_T\right)h_{uT}. \qedhere
\end{align*}
\end{proof}

\begin{corollary}[Nonvanishing terminal difference]
\label{cor:zero-b-nonzero-h}
Under \cref{ass:top-degree-equality}, $h\ne0$.
\end{corollary}
\begin{proof}
If $h=0$, \cref{lem:interface-identities}\ref{it:interface-master} gives
$(L+1)b_T=0$ for every $T$, hence $b=0$.
The interpolation formulas of
\cref{lem:interface-identities}\ref{it:interface-terminal}--\ref{it:interface-no-terminal}
then give $f=0$ on every $k$-subset, a contradiction.
\end{proof}

We now reverse the interface identities.  The following reconstruction is purely
algebraic: compatible data $(h,b,p)$ uniquely determine a function in
$\cH_k$.  The additional conditions for producing an equality vector are recorded in \cref{cor:reconstruction-equality}.

\begin{lemma}[Algebraic reconstruction from interface data]\label{lem:interface-reconstruction}
Fix $3\le k\le\lfloor \frac{n}{2}\rfloor$, put $W=V\setminus\{x,y\}$,
$N=|W|$, and $L=N-k+2$.  Fix a function
$p:W\to[0,1]$, and write $p_u:=p(u)$ and $S_A:=\sum_{u\in A}p_u$.
Suppose $h:\binomset W{k-1}\to\RR$ and $b:\binomset W{k-2}\to\RR$
satisfy \cref{lem:interface-identities} items \ref{it:interface-down} (the down-map identities), \ref{it:interface-symmetry} (symmetry), \ref{it:interface-consistency} (consistency), and \ref{it:interface-master} (master identity).

For $T\in\binomset W{k-2}$, $\xi\in\binomset W{k-1}$,
$u\in\xi$, $Q\in\binomset Wk$, and $z\in Q$, define $f: \binomset V{k} \to \mathbb{R}$ by
\begin{align}
  f(xyT)&=b_T,\label{eq:reconstruction-both}\\
  f(x\xi)&=b_{\xi\setminus\{u\}}+(1-p_u)h_\xi,\label{eq:reconstruction-x}\\
  f(y\xi)&=b_{\xi\setminus\{u\}}-p_uh_\xi,\label{eq:reconstruction-y}\\
  f(Q)&=(1-p_z)f\bigl(x(Q\setminus\{z\})\bigr)
          +p_zf\bigl(y(Q\setminus\{z\})\bigr).\label{eq:reconstruction-interior}
\end{align}
Then:
\begin{enumerate}[label=\textup{(\roman*)},ref=(\roman*)]
\item\label{it:reconstruction-well-defined}
Equations \eqref{eq:reconstruction-both}--\eqref{eq:reconstruction-interior}
are independent of the choices of $u$ and $z$, and uniquely define
$f$.

\item\label{it:reconstruction-top}
The reconstructed function $f$ belongs to $\cH_k$.

\item\label{it:reconstruction-norm}
The terminal norm of $f$ satisfies
$\delta_f(x,y)^2=\binom nk^{-1}\sum_{\xi\in\binomset W{k-1}}h_\xi^2$.

\item\label{it:reconstruction-distances}
For every $z\in W$,
$\delta_f(x,z)=p_z\delta_f(x,y)$ and
$\delta_f(z,y)=(1-p_z)\delta_f(x,y)$.
\end{enumerate}
\end{lemma}

\begin{proof}
Item~\ref{it:reconstruction-well-defined}: First let
$\xi\in\binomset W{k-1}$ and choose distinct $u,v\in\xi$.
With $R=\xi\setminus\{u,v\}$,
\cref{lem:interface-identities}\ref{it:interface-symmetry} gives
$b_{vR}-b_{uR}=(p_u-p_v)h_\xi$.  Consequently,
\[
\begin{aligned}
  b_{vR}+(1-p_u)h_\xi
    &=b_{uR}+(1-p_v)h_\xi,\\
  b_{vR}-p_uh_\xi
    &=b_{uR}-p_vh_\xi.
\end{aligned}
\]
Thus \eqref{eq:reconstruction-x} and \eqref{eq:reconstruction-y}
are independent of the chosen vertex of $\xi$.

Next let $Q\in\binomset Wk$, choose distinct $z,w\in Q$, and put
$T=Q\setminus\{z,w\}$.  Evaluating \eqref{eq:reconstruction-interior} with $z$ and then with $w$
gives, respectively,
\[
\begin{aligned}
 (1-p_z)f(xwT)+p_zf(ywT)
   &=(1-p_z)\bigl[b_T+(1-p_w)h_{wT}\bigr]
     +p_z\bigl[b_T-p_wh_{wT}\bigr]\\
   &=b_T+(1-p_z-p_w)h_{wT},\\
 (1-p_w)f(xzT)+p_wf(yzT)
   &=b_T+(1-p_z-p_w)h_{zT}.
\end{aligned}
\]
\cref{lem:interface-identities}\ref{it:interface-consistency} says
$(1-p_z-p_w)(h_{zT}-h_{wT})=0$, so these evaluations agree.
Every $k$-subset contains both terminals, exactly one terminal, or
neither terminal.  Therefore under the hypotheses on the interface data $h$ and $b$, \eqref{eq:reconstruction-both}--\eqref{eq:reconstruction-interior}
defines $f$ unambiguously and uniquely.

Item~\ref{it:reconstruction-top}: Fix $T\in\binomset W{k-2}$.
For each $z\in W\setminus T$, 
\cref{lem:interface-identities}\ref{it:interface-master} is equivalent to
\begin{align}
\label{eq:hzTeq}
  (S_W-S_T-L+Lp_z)h_{zT}=(L+1)b_T.
\end{align}
Evaluating the down-map identity $D^W_{k-1}h=0$ from
\cref{lem:interface-identities}\ref{it:interface-down} at the $(k-2)$-subset $T$ gives
$\sum_{z\in W\setminus T}h_{zT}=0$.
Using this identity and the fact that $|W\setminus T|=L$, we sum the equation \eqref{eq:hzTeq} over all $z\in W\setminus T$ to obtain
\[
  \sum_{z\in W\setminus T}p_zh_{zT}=(L+1)b_T.
\]
Equations \eqref{eq:reconstruction-both}--\eqref{eq:reconstruction-y}
now give, for the $(k-1)$-subsets $xT$ and $yT$,
\[
\begin{aligned}
 (D_kf)(xT)
   &=f(xyT)+\sum_{z\in W\setminus T}f(xzT)
     =(L+1)b_T+\sum_{z\in W\setminus T}(1-p_z)h_{zT}=0,\\
 (D_kf)(yT)
   &=f(xyT)+\sum_{z\in W\setminus T}f(yzT)
     =(L+1)b_T-\sum_{z\in W\setminus T}p_zh_{zT}=0.
\end{aligned}
\]
For $R\in\binomset W{k-3}$, the remaining case containing
both terminals follows from the other down-map identity $D^W_{k-2}b=0$ in
\cref{lem:interface-identities}\ref{it:interface-down}:
\[
 (D_kf)(xyR)=\sum_{z\in W\setminus R}b_{zR}
              =(D_{k-2}^Wb)(R)=0.
\]
Finally, let $\xi\in\binomset W{k-1}$, choose $u\in\xi$, and put
$T=\xi\setminus\{u\}$.  For every $z\in W\setminus\xi$,
\eqref{eq:reconstruction-x}--\eqref{eq:reconstruction-interior} give
$f(z\xi)=b_T+(1-p_u-p_z)h_\xi$.
As $|W\setminus\xi|=L-1$, it follows that
\[
 (D_kf)(\xi)
  =f(x\xi)+f(y\xi)+\sum_{z\in W\setminus\xi}f(z\xi)
  =(L+1)b_T+(L-S_W-Lp_u+S_T)h_\xi=0,
\]
where the final equality is exactly
\cref{lem:interface-identities}\ref{it:interface-master}.
These four types exhaust the $(k-1)$-subsets of $V$.
Hence $D_kf=0$, or $f\in\cH_k$.

Item~\ref{it:reconstruction-norm}: Subtracting \eqref{eq:reconstruction-y} from \eqref{eq:reconstruction-x}
gives $f(x\xi)-f(y\xi)=h_\xi$.  Applying the coordinate formula
\eqref{eq:delta-coordinate} to the pair $x,y$ yields the claimed formula.

Item~\ref{it:reconstruction-distances}: Fix $z\in W$.
In \eqref{eq:delta-coordinate} for the pair $x,z$, the indexing
$(k-1)$-subsets of $V\setminus\{x,z\}$ either avoid $y$ or contain $y$.
Writing them, respectively, as
$\xi\in\binomset{W\setminus\{z\}}{k-1}$ and
$yR$ where $R\in\binomset{W\setminus\{z\}}{k-2}$, gives
\[
 \delta_f(x,z)^2=\binom nk^{-1}\left[
 \sum_{\xi\in\binomset{W\setminus\{z\}}{k-1}}
        [f(x\xi)-f(z\xi)]^2
 +\sum_{R\in\binomset{W\setminus\{z\}}{k-2}}
        [f(xyR)-f(zyR)]^2\right].
\]
Equations \eqref{eq:reconstruction-x}--\eqref{eq:reconstruction-interior}  give
\[
 f(x\xi)-f(z\xi)=p_zh_\xi,
 \qquad
 f(xyR)-f(zyR)=p_zh_{zR}.
\]
Consequently,
\[
 \delta_f(x,z)^2
   =p_z^2\binom nk^{-1}\left[
       \sum_{\xi\in\binomset{W\setminus\{z\}}{k-1}}h_\xi^2
       +\sum_{R\in\binomset{W\setminus\{z\}}{k-2}}h_{zR}^2
     \right]
   =p_z^2\delta_f(x,y)^2.
\]
The second equality uses the partition of the $(k-1)$-subsets of
$W$ into those avoiding $z$ and those containing $z$, together
with item~\ref{it:reconstruction-norm}.

For the pair $z,y$, split the indexing subsets in
\eqref{eq:delta-coordinate} according to whether they contain $x$.
The corresponding differences are
\[
 f(z\xi)-f(y\xi)=(1-p_z)h_\xi,
 \qquad
 f(zxR)-f(yxR)=(1-p_z)h_{zR}.
\]
The same partition of the $(k-1)$-subsets of $W$ therefore gives
$\delta_f(z,y)^2=(1-p_z)^2\delta_f(x,y)^2$.
Since $0\le p_z\le1$, taking square roots proves
$\delta_f(x,z)=p_z\delta_f(x,y)$ and
$\delta_f(z,y)=(1-p_z)\delta_f(x,y)$.
\end{proof}

\begin{corollary}[Equality criterion for reconstructed functions]
\label{cor:reconstruction-equality}
Let $f\in\cH_k$ be reconstructed from compatible interface data
$(h,b,p)$ as in \cref{lem:interface-reconstruction}.  Then:
\begin{enumerate}[label=\textup{(\roman*)},ref=(\roman*)]
\item \label{it:f=0-iff-h=0} $f=0$ if and only if $h=0$.
\item \label{it:Kf=0} Suppose $h\ne0$ and $p$ is the harmonic $x$--$y$ voltage of $G$,
restricted to $W$.  Then
\[
  \cK_{xy}^{(k)}f=0
  \quad\Longleftrightarrow\quad
  \delta_f(u,v)=\delta_f(x,y)|p_u-p_v|
  \quad\text{for every }\{u,v\}\subset W\text{ with }c_{uv}>0.
\]
\end{enumerate}
\end{corollary}

\begin{proof}
Item~\ref{it:f=0-iff-h=0}: If $h=0$, the master identity
\cref{lem:interface-identities}\ref{it:interface-master} gives $b=0$,
and \eqref{eq:reconstruction-both}--\eqref{eq:reconstruction-interior} then give $f=0$.  
Conversely,
\eqref{eq:reconstruction-x}--\eqref{eq:reconstruction-y} gives
$h_\xi=f(x\xi)-f(y\xi)$, so $f=0$ implies $h=0$.

Item~\ref{it:Kf=0}: Suppose $h\ne0$.  By
\cref{lem:interface-reconstruction}\ref{it:reconstruction-norm},
$\delta_f(x,y)>0$.  The distance identities in
\cref{lem:interface-reconstruction}\ref{it:reconstruction-distances}
then give $\phi_f=p$ on $V$, with $p(x)=0$ and $p(y)=1$.
Since $p$ is harmonic, \cref{prop:MPL-equality,cor:operator-equality}
reduce residual-kernel membership to the positive-edge distance
identities.  Those identities already hold on the edges incident to either
terminal by \cref{lem:interface-reconstruction}\ref{it:reconstruction-distances};
and also on $\{x,y\}$ if it is a positive edge. 
Therefore it remains to check them on the positive edges within $W$, which is exactly the right-hand side criterion in the two-way implication.
\end{proof}

\subsection{Voltage rigidity} \label{sec:voltage-rigidity}

We show that a nonzero equality vector in $\cH_k$ permits at most two
nonterminal voltage values.  For $k=3$, eliminating the interface variables
yields a determinant; for $k\ge4$, comparing the master identity across four
vertices produces a Vandermonde factor in three voltage variables.
Retain $L=N-k+2$ from \cref{lem:interface-identities}\ref{it:interface-master}.

\subsubsection{Degree three}

When $k=3$, the interface data \eqref{eq:interface-data} are $h_{uv} := f(xuv)-f(yuv)$ and $b_u:= f(xyu)$, where $u,v\in W$ are distinct.  Then $L=N-1$, and
\cref{lem:interface-identities}\ref{it:interface-master} becomes
\begin{equation}\label{eq:k3-master}
  N b_v+
  \bigl(\Theta+p_v-(N-1)p_u\bigr)h_{uv}=0,
  \qquad
  \Theta:=N-1-S_W.
\end{equation}
By symmetry, $Nb_u + (\Theta+p_u - (N-1)p_v) h_{uv}=0$.
Hence $b_u$ and $b_v$ satisfy the proportionality relation
\begin{equation}\label{eq:k3-b-relation}
  \bigl(\Theta+p_u-(N-1)p_v\bigr)b_v
  =\bigl(\Theta+p_v-(N-1)p_u\bigr)b_u.
\end{equation}

\begin{lemma}[Degree-$3$ determinant]\label{lem:k3-determinant}
For three distinct vertices $u,v,w\in W$, the determinant of the three pair relations
\eqref{eq:k3-b-relation} equals
\begin{equation}\label{eq:k3-determinant}
  -N(N-1)(p_u-p_v)(p_u-p_w)(p_v-p_w).
\end{equation}
Hence, if $p_u,p_v,p_w$ are pairwise distinct, then
$b_u=b_v=b_w=0$.
\end{lemma}

\begin{proof}
Writing the relations for $(u,v)$, $(u,w)$, and $(v,w)$ as a homogeneous system in
$(b_u,b_v,b_w)$ gives the coefficient matrix
\[
M=\begin{pmatrix}
 -(\Theta+p_v-(N-1)p_u) & \Theta+p_u-(N-1)p_v & 0\\
 -(\Theta+p_w-(N-1)p_u) & 0 & \Theta+p_u-(N-1)p_w\\
 0 & -(\Theta+p_w-(N-1)p_v) & \Theta+p_v-(N-1)p_w
\end{pmatrix}.
\]
Treat $\Theta$ temporarily as an
independent scalar and put $q=N-1$.  Expanding the determinant using, say, the cofactor formula gives
\[
\begin{aligned}
\det M={}&(\Theta+p_u-qp_v)(\Theta+p_w-qp_u)(\Theta+p_v-qp_w)\\
 &-(\Theta+p_v-qp_u)(\Theta+p_u-qp_w)(\Theta+p_w-qp_v).
\end{aligned}
\]
Interchanging any two of $p_u,p_v,p_w$ exchanges these products and hence
changes the sign of $\det M$.  Since $\det M$ has total degree at most three in the three
voltages, it is a constant multiple of
$(p_u-p_v)(p_u-p_w)(p_v-p_w)$.  The coefficient of $p_u^2p_v$
is $-q$ in the first product and $q^2$ in the second, so its coefficient
in $\det M$ is $-q-q^2=-N(N-1)$.  This proves \eqref{eq:k3-determinant}
for every $\Theta$, in particular for $\Theta=N-1-S_W$.
\end{proof}

\subsubsection{Degrees at least four}

\begin{lemma}[Four-point voltage identity]\label{lem:four-point}
Assume $k\ge4$.  Let $R\in\binomset W{k-4}$, and $u,v,w,z$ be four
distinct vertices of $W\setminus R$.  Then
$
  L(L+1)
  (p_v-p_w)(p_v-p_z)(p_w-p_z)
  h_{uvwR}=0
  $.
\end{lemma}

\begin{proof}
For distinct $a,b,c\in W\setminus R$, set
$
  K_a(b,c):=L-S_W+S_R-Lp_a+p_b+p_c
 $.
Apply \cref{lem:interface-identities}\ref{it:interface-master} with $T=uvR$,
first choosing the vertex $w\in W\setminus T$, and then $z\in W\setminus T$.
The two equations contain the same term $(L+1)b_{uvR}$, which cancels
upon subtraction and gives
\[
K_w(u,v)h_{uvwR}=K_z(u,v)h_{uvzR}.
\]
Likewise, taking $T=uzR$ with the respective choices $v,w\in W\setminus T$,
and then $T=uwR$ with $z,v\in W\setminus T$, gives
\begin{align*}
  K_v(u,z)h_{uvzR}&=K_w(u,z)h_{uwzR},\\
  K_z(u,w)h_{uwzR}&=K_v(u,w)h_{uvwR}.
\end{align*}
Multiplying these identities so as to eliminate $h_{uvzR}$ and $h_{uwzR}$, without
dividing by any coefficient, yields
\[
  0=\bigl[
  K_v(u,z)K_w(u,v)K_z(u,w)
 -K_w(u,z)K_v(u,w)K_z(u,v)
  \bigr]h_{uvwR}.
\]
The bracket is a polynomial of total degree at most three in $p_v,p_w,p_z$ and
changes sign when any two of these three voltages are exchanged.  It therefore
vanishes whenever two of them coincide and is a constant multiple of
$(p_v-p_w)(p_v-p_z)(p_w-p_z)$.  The coefficient of $p_v^2p_w$ is $-L(L+1)$, so the
bracket equals
$-L(L+1)(p_v-p_w)(p_v-p_z)(p_w-p_z)$.  This proves the lemma.
\end{proof}

\begin{proposition}[At most two nonterminal voltage values]
\label{prop:two-voltage-rigidity}
Under \cref{ass:top-degree-equality}, the set
$\{p_u:u\in W\}$ has cardinality at most two.
\end{proposition}

\begin{proof}
Suppose, on the contrary, that $\{p_u: u\in W\}$ contains at least three distinct voltage values.
We discuss the two alternatives, $k=3$ and $k\ge 4$.

\underline{The case $k=3$:}  Fix $u\in W$.  Since at least three voltage values
occur, we can choose $v,w\in W$ such that $p_u,p_v,p_w$ are pairwise
distinct. 
\cref{lem:k3-determinant} gives
$b_u=b_v=b_w=0$.  As $u$ was arbitrary, $b_u=0$ for every $u\in W$.
Hence $b=0$, and $h\ne0$ by \cref{cor:zero-b-nonzero-h}.  
By \cref{lem:interface-identities}\ref{it:interface-symmetry},
$(p_u-p_v)h_{uv}=0$ for distinct $u,v\in W$.  In particular,
$h_{uv}=0$ whenever $p_u\ne p_v$.

For distinct $u,v$ with $p_u=p_v=t$, we use $b_v=0$ to reduce equation \eqref{eq:k3-master} to $(\Theta-(N-2)t)h_{uv}=0$.  The quantity
$\Theta=N-1-S_W$ is independent of $u,v$, and the standing bound
$k=3\le \frac{n}{2}$ gives $N=n-2\ge4$.  Hence $N-2>0$, so $h_{uv}\ne0$
for a pair with equal voltages is possible only at the common value
$
  t_*:=\frac{\Theta}{N-2}$.
Thus $h_{uv}$ can be nonzero only when $p_u=p_v=t_*$.

\underline{The case $k\ge 4$:}  We first show that $h_\xi=0$ whenever
$\xi\in\binomset W{k-1}$ contains vertices with different voltages.
If exactly two voltage values occur on $\xi$, choose $d\in W\setminus\xi$
at a third voltage value, and choose $b,c\in\xi$ at different voltages.
If at least three values occur on $\xi$, choose any $d\in W\setminus\xi$
and then $b,c\in\xi$ so that $p_b,p_c,p_d$ are pairwise distinct.
These choices are possible: in the first case the third value occurs
outside $\xi$, and in the second case at least two values on $\xi$
differ from $p_d$. 
Moreover, $W\setminus\xi\ne\varnothing$ because
$N\ge2k-2>k-1$.
Choose $a\in\xi\setminus\{b,c\}$ and put
$R=\xi\setminus\{a,b,c\}$.  Applying \cref{lem:four-point} with
$(u,v,w,z)=(a,b,c,d)$ gives
$
 L(L+1)(p_b-p_c)(p_b-p_d)(p_c-p_d)h_\xi=0$.
Since $L=N-k+2\ge k\ge4$ and $p_b,p_c,p_d$ are pairwise distinct,
the scalar coefficient is nonzero.  Thus $h_\xi=0$, as claimed.

Next fix $T\in\binomset W{k-2}$.  Choose $a\in W\setminus T$ so that
$p$ is not constant on $aT$: any $a\in W\setminus T$ works if $p$
is already nonconstant on $T$, while otherwise we choose $a$ at a
voltage different from the common value on $T$.
The preceding paragraph gives $h_{aT}=0$, and
\cref{lem:interface-identities}\ref{it:interface-master} gives
$(L+1)b_T=0$.  Therefore $b_T=0$ for every $T$, or $b=0$.
Hence $h\neq 0$ by \cref{cor:zero-b-nonzero-h}.

It remains to consider $\xi\in\binomset W{k-1}$ on which $p$ is
constant with common voltage value $t$.  Taking any $u\in\xi$ and
$T=\xi\setminus\{u\}$ in \cref{lem:interface-identities}\ref{it:interface-master} now yields $
 (L-S_W+(2k-N-4)t)h_\xi=0$.
Since $N\ge2k-2$, we have $N+4-2k\ge2$.  Thus a nonzero $h_\xi$
on such a set requires its common voltage to equal
$t_*:=\frac{L-S_W}{N+4-2k}$.  
Thus, also in this case, $h_\xi$ can be nonzero only if
$p$ is constant with value $t_*$ on $\xi$.

\emph{Common contradiction.}
In both cases we have established $b=0$, $h\ne0$, and a scalar $t_*$
(depending on the case) such that $h_\xi\ne0$ only if every vertex of
$\xi$ has voltage $t_*$.  Choose such a $\xi$ with $h_\xi\ne0$, fix
$u\in\xi$, and put $T=\xi\setminus\{u\}$.  For any $d\in W$ with
$p_d\ne t_*$, we have $d\notin\xi$.  The set $dT$ contains vertices
with two distinct voltage values, so the support restriction proved in
both cases gives $h_{dT}=0$.  As $b=0$,
\cref{lem:interface-identities}\ref{it:interface-terminal} gives
\[
 f(xdT)=f(ydT)=0,\qquad
 f(x\xi)=(1-t_*)h_\xi,\qquad f(y\xi)=-t_*h_\xi.
\]
Applying \cref{lem:interface-identities}\ref{it:interface-no-terminal}
with distinguished vertex $u$ and remaining set $dT$ gives
$f(d\xi)=0$.  Applying the same identity again
with distinguished vertex $d$ and remaining set $\xi$ gives
\[
 0=f(d\xi)=(1-p_d)f(x\xi)+p_df(y\xi)
          =(1-p_d-t_*)h_\xi.
\]
Since $h_\xi\ne0$, every voltage different from $t_*$ must equal
$1-t_*$.  This contradicts the assumption of at least three distinct
nonterminal voltage values.
\end{proof}

\section{Quadratic equality kernel}
\label{sec:quadratic-classification}

\cref{prop:two-voltage-rigidity} established that a nonzero equality
vector of Hoeffding degree $k\ge3$ permits at most two nonterminal voltage
values.  
In this section we analyze nonzero equality vectors of Hoeffding degree $k=2$. 
The analysis produces a sharper classification: either every
nonterminal voltage equals $\frac12$, or exactly two complementary values $a$
and $1-a$ occur.

We use the following notation throughout this section.
If $n\le3$, then
$\lfloor \frac{n}{2}\rfloor\le1$, so there is no Hoeffding degree-two layer to classify.
Hence throughout this section assume $n=|V|\ge4$.  
As before, set $W:=V\setminus\{x,y\}$ and $N:=|W|=n-2$.
Let $G[W]$ be the graph on $W$ obtained by retaining exactly the edges of positive
conductance, and let $W_1,\dots,W_q$ be the connected components of $G[W]$.

\subsection{Terminal balance and quadratic interface identities}

Recall the harmonic voltage condition \eqref{eq:harmonic-voltage-equation}. 
If the harmonic voltage is constant, say equal to $\rho_i$, on a component $W_i$,
then the internal-edge contributions cancel in the harmonicity equation.  Harmonicity
at $u\in W_i$ therefore reduces to
$\rho_i c_{xu}=(1-\rho_i)c_{yu}$.  This motivates the following definition.

\begin{definition}[Terminal balance]
A component $W_i$ is \emph{terminal-balanced} at voltage $\rho_i\in[0,1]$ if
\begin{equation}\label{eq:terminal-balance}
  \rho_i c_{xu}=(1-\rho_i)c_{yu}
  \qquad\text{for every }u\in W_i.
\end{equation}
\end{definition}

Because $G$ is connected, every $W_i$ has at least one vertex incident to $x$ or $y$.
Hence if \eqref{eq:terminal-balance} holds, $\rho_i$ is unique.  Equivalently, at every vertex with $c_{xu}+c_{yu}>0$,
$\rho_i=c_{yu}/(c_{xu}+c_{yu})$, and this ratio is constant over all vertices of $W_i$ adjacent to a terminal.

The following lemma records the identities satisfied by every quadratic equality
vector.
Let $F\in\cH_2\setminus\{0\}$ satisfy $\cK_{xy}^{(2)}F=0$.  Since a nonzero
constant does not belong to $\cH_2$, \cref{cor:operator-equality} gives
$\delta_F(x,y)>0$.  Let $p_u:=\phi_F(u)$ be the harmonic voltage from
\cref{prop:MPL-equality}, and following the notation of \cref{sec:interface}, set $S_W:=\sum_{u\in W}p_u$. 
The interface functions from \eqref{eq:interface-data} specialize to $h_u:=F(xu)-F(yu)$ and  $b:=F(xy)$.

\begin{lemma}[Scalar quadratic interface]\label{lem:q-scalar-interface}
 The following hold for distinct $u,v\in W$:
\begin{enumerate}[label=\textup{(\roman*)},ref=(\roman*)]
\item\label{it:q-interpolation}
The interpolation identities are
$F(xu)=b+(1-p_u)h_u$, $F(yu)=b-p_uh_u$, and
$F(uv)=b+(1-p_u-p_v)h_v=b+(1-p_u-p_v)h_u$.

\item\label{it:q-consistency}
The two-vertex consistency identity is $(1-p_u-p_v)(h_u-h_v)=0$.

\item\label{it:q-zero-sum}
The interface vector has zero sum: $\sum_{u\in W}h_u=0$.

\item\label{it:q-master}
The scalar master identity is $(N+1)b+(N-S_W-Np_u)h_u=0$.

\item\label{it:q-norm}
The terminal norm squared is
$\delta_F(x,y)^2=\binom n2^{-1}\sum_{u\in W}h_u^2$. In particular $h\neq0$.
\end{enumerate}
\end{lemma}

\begin{proof}
By \cref{lem:every-vertex-geodesic}, every $u\in W$ satisfies
$\delta_F(x,y)=\delta_F(x,u)+\delta_F(u,y)$.  Applying
\eqref{eq:terminal-coordinate-proportionality} to the corresponding coordinates of
$\mathbf d_F^{xy}$ gives the first two formulas in item~\ref{it:q-interpolation}.
Applying it to the coordinate indexed by $v$ gives the first expression for $F(uv)$;
interchanging $u$ and $v$ gives the second.  Their difference proves
item~\ref{it:q-consistency}.

The relations $D_2F(x)=D_2F(y)=0$ say
$b+\sum_{u\in W}F(xu)=0$ and $b+\sum_{u\in W}F(yu)=0$; subtracting and
using item~\ref{it:q-interpolation} gives $\sum_{u\in W}h_u=0$.
For fixed $u\in W$, the remaining down-map identity, together with the
interpolation formulas, gives
\[
\begin{aligned}
0&=F(xu)+F(yu)+\sum_{v\in W\setminus\{u\}}F(uv)\\
 &= (N+1)b+
 \left[1-2p_u+\sum_{v\in W\setminus\{u\}}(1-p_u-p_v)\right]h_u\\
 &= (N+1)b+(N-S_W-Np_u)h_u,
\end{aligned}
\]
where the last step uses
$\sum_{v\in W\setminus\{u\}}p_v=S_W-p_u$.
This proves item~\ref{it:q-master}.  Finally, item~\ref{it:q-norm}
follows from \eqref{eq:delta-coordinate} for the pair $x,y$.
Since $F\in\cH_2$ is a nonzero equality vector,
$\delta_F(x,y)>0$ by \cref{cor:operator-equality}, so $h\neq0$.
\end{proof}

We distinguish the \emph{midpoint case} $b=0$ from the
\emph{complementary case} $b\ne0$.  These names describe the harmonic
voltage functions that will be proved below: the first case forces $p\equiv\frac12$ on $W$, while the second case forces
exactly two values $a$ and $1-a$, with $0\le a<\frac12$.

\subsection{The midpoint case: \texorpdfstring{$b=0$}{b=0}}

\begin{lemma}[Midpoint case]\label{lem:q-midpoint-alternative}
If $b=0$, then $p_u=\frac12$ for every $u\in W$,
and there is a function $\alpha:W\to\RR$ with $\sum_{u\in W}\alpha(u)=0$ such that for all $u,v \in W$ with $u\neq v$,
\begin{equation}\label{eq:q-midpoint-form}
  F(xu)=\alpha(u),
  \qquad
  F(yu)=-\alpha(u),
  \qquad
  F(xy)=0,\qquad F(uv)=0.
\end{equation}
Moreover, $\alpha$ is constant on each component $W_i$.
\end{lemma}

\begin{proof}
When $b=0$, \cref{lem:q-scalar-interface}\ref{it:q-master} shows that
$h_u\neq0$ requires $p_u=\frac{1}{N}(N-S_W)=:t$.
Since $h\neq0$ and $\sum_{w\in W}h_w=0$ by
\cref{lem:q-scalar-interface}\ref{it:q-zero-sum}, there are vertices
$u,v\in W$ with $h_u>0$ and $h_v<0$, and therefore $p_u=p_v=t$.
The consistency identity \cref{lem:q-scalar-interface}\ref{it:q-consistency}
then becomes
$
  (1-2t)(h_u-h_v)=0$.
As $h_u-h_v>0$, we obtain $t=\frac12$.
If $h_w=0$, applying the same consistency identity to $u,w$ gives
$(1-p_u-p_w)h_u=0$, so $p_w=1-p_u=\frac12$ as well.

Set $\alpha(u):=\frac{h_u}{2}$.  Using \cref{lem:q-scalar-interface}\ref{it:q-interpolation} gives
\eqref{eq:q-midpoint-form}, and \cref{lem:q-scalar-interface}\ref{it:q-zero-sum} gives
$\sum_u\alpha(u)=0$.  For $u,v\in W$, direct use of
\eqref{eq:delta-coordinate} gives
\begin{equation}\label{eq:q-midpoint-W-distance}
  \delta_F(u,v)^2
  =2\binom n2^{-1}\bigl(\alpha(u)-\alpha(v)\bigr)^2.
\end{equation}
If $c_{uv}>0$, then \eqref{eq:edge-equality} and $p_u=p_v=\frac12$ give
$\delta_F(u,v)=0$.  Hence $\alpha(u)=\alpha(v)$ on every positive edge inside $W$,
and therefore on every component $W_i$.
\end{proof}

\subsection{The complementary case: \texorpdfstring{$b\neq0$}{b not equal to 0}}

We first construct a quadratic equality vector when the two nonempty
voltage classes $A$ and $B$ carry values $a$ and $1-a$ and no positive edge
joins them.  The construction assigns one value to each pair type determined
by $x,y,A$, and $B$; we then establish necessity and uniqueness up to a scalar multiple.

\begin{lemma}[Quadratic construction in the complementary case]
\label{lem:complementary-quadratic-construction}
Suppose $W=A\sqcup B$ with $A,B$ nonempty, and let $0\le a<\frac12$.
Assume that the harmonic $x$--$y$ voltage takes value $a$ on $A$ and $1-a$ on $B$, and that
no positive edge joins $A$ to $B$.  Then there exists a nonzero $F\in\cH_2$ satisfying
$\cK_{xy}^{(2)}F=0$.
\end{lemma}

\begin{proof}
Fix a scalar $\Delta\neq0$ and put $R_0:=\frac{|A|(2a-1)}{N+1}\Delta$.
Define $F:\binomset{V}{2}\to\RR$ by the pair-value table
\eqref{eq:complementary-pair-table} below.
To explain the notation, $F(xA)$ denotes
the common value $F(xu)$ for $u\in A$, and $F(AB)$ denotes the common
value $F(uv)$ for $u\in A$, $v\in B$.  The $AA$ row is used only if $|A|\ge2$, and the $BB$
row only if $|B|\ge2$.
\begin{equation}\label{eq:complementary-pair-table}
\begin{array}{c|c}
\text{pair type for } zw&F(zw)\\ \hline
xy,\ AB&R_0\\[1mm]
xA&R_0+(1-a)\Delta\\
yA&R_0-a\Delta\\
xB&R_0-\dfrac{a |A|}{|B|}\Delta\\
yB&R_0+\dfrac{(1-a)|A|}{|B|}\Delta\\
AA&R_0+(1-2a)\Delta\\
BB&R_0+\dfrac{|A|}{|B|}(1-2a)\Delta
\end{array}
\end{equation}

The rest of the proof verifies that $F$ is a nonzero quadratic equality vector.
To verify $D_2F=0$, first compute the row indexed by $x$:
\[
 D_2F(x)
 =F(xy)+|A|F(xA)+|B|F(xB)
 =(N+1)R_0+|A|(1-2a)\Delta=0.
\]
For $u\in A$, the $|A|-1$ other vertices of $A$ and all $|B|$
vertices of $B$ contribute
\[
 D_2F(u)
 =F(xA)+F(yA)+\sum_{w\in A\setminus\{u\}}F(uw)+|B|F(AB)
 =(N+1)R_0+|A|(1-2a)\Delta=0.
\]
The sum over $A\setminus\{u\}$ is empty when $|A|=1$.  For $v\in B$,
the remaining two kinds of row sums have the forms
\[
\begin{aligned}
 D_2F(y)&=F(xy)+|A|F(yA)+|B|F(yB),\\
 D_2F(v)&=F(xB)+F(yB)+|A|F(AB)
            +\sum_{w\in B\setminus\{v\}}F(vw).
\end{aligned}
\]
Substituting the pair values from \eqref{eq:complementary-pair-table} in both rows gives
\[
  D_2F(y)=D_2F(v)=(N+1)R_0+|A|(1-2a)\Delta=0.
\]
Hence $D_2F=0$.

Set $\Lambda:=\delta_F(x,y)^2$.  The coordinates of the difference
for $(x,y)$ are indexed by $W$, so \eqref{eq:delta-coordinate} gives
\[
 \binom n2\Lambda
 =\left(|A|+\frac{|A|^2}{|B|}\right)\Delta^2
 =\frac{N|A|}{|B|}\Delta^2.
\]
For $u\in A$, the coordinates indexed by $y$ and by
$A\setminus\{u\}$ contribute $|A|$ differences of magnitude
$a|\Delta|$, while those indexed by $B$ contribute $|B|$
differences of magnitude $\frac{a|A|}{|B|}|\Delta|$.
For $v\in B$, the analogous differences for $(x,v)$ have magnitudes
$(1-a)|\Delta|$ on $A$ and $\frac{(1-a)|A|}{|B|}|\Delta|$ on
$\{y\}\cup(B\setminus\{v\})$.  Consequently,
\[
\begin{aligned}
\binom n2\delta_F(x,u)^2
 &=a^2\left(|A|+\frac{|A|^2}{|B|}\right)\Delta^2
 =a^2\binom n2\Lambda,\\
\binom n2\delta_F(x,v)^2
 &=(1-a)^2\left(|A|+\frac{|A|^2}{|B|}\right)\Delta^2
 =(1-a)^2\binom n2\Lambda.
\end{aligned}
\]
The same calculation with $y$ instead of $x$ gives
$\delta_F(y,u)^2=(1-a)^2\Lambda$ and
$\delta_F(y,v)^2=a^2\Lambda$.
The distances within either class vanish because the corresponding
coordinate differences are zero for every pair type.

For $u\in A$ and $v\in B$, the coordinates indexed by $x$, $y$, the
$|A|-1$ other vertices of $A$, and the $|B|-1$ other vertices of $B$ give,
respectively, the four squared differences below:
\[
\begin{aligned}
  \binom n2\frac{\delta_F(u,v)^2}{\Delta^2}
  &=\left(1-a+a\frac{|A|}{|B|}\right)^2
    +\left(a+(1-a)\frac{|A|}{|B|}\right)^2
  +(|A|-1)(1-2a)^2
    +(|B|-1)\left(\frac{|A|}{|B|}\right)^2(1-2a)^2\\
  &=(1-2a)^2\frac{N|A|}{|B|}
       +\frac{2a(1-a)N^2}{|B|^2}.
\end{aligned}
\]
Consequently,
\begin{equation}\label{eq:F-dist-AB}
  \delta_F(u,v)^2
  =(1-2a)^2\Lambda
   +\binom n2^{-1}\frac{2a(1-a)N^2}{|B|^2}\Delta^2.
\end{equation}
By hypothesis no positive edge joins $A$ to $B$.  The displayed distances give
\eqref{eq:edge-equality} on every positive edge joining a terminal to $W$, and the
zero within-block distances give it on internal positive edges.  On $\{x,y\}$, if
it is an edge, the equality is tautological.  Moreover, the normalized distance from
$x$ agrees with the harmonic voltage $p$.  For any $u\in A$, $F(xu)-F(yu)=\Delta\ne0$, so
$\delta_F(x,y)>0$.  Since $\phi_F=p$ and the positive-edge conditions
have been verified, \cref{prop:MPL-equality,cor:operator-equality} give
$0\ne F\in\ker(\cK_{xy}^{(2)}|_{\cH_2})$.
\end{proof}

\begin{lemma}[Complementary case]\label{lem:q-complementary-alternative}
Suppose $b\neq0$.  Then:
\begin{enumerate}[label=\textup{(\roman*)},ref=(\roman*)]
\item\label{it:q-complementary-voltages}
There are exactly two nonterminal voltage values $a$ and $1-a$, where
$0\le a<\frac12$.

\item\label{it:q-complementary-normal-form}
Writing $A:=\{u\in W:p_u=a\}$ and $B:=\{u\in W:p_u=1-a\}$, the space of quadratic equality vectors with this voltage
partition is one-dimensional.  More precisely, if $\Delta:=h_u$ on $A$, then
$h=-\frac{|A|}{|B|}\Delta$ on $B$, and $F$ is exactly the function in
\eqref{eq:complementary-pair-table} with
$R_0=\frac{|A|(2a-1)}{N+1}\Delta$.

\item\label{it:q-complementary-separation}
No positive edge inside $W$ joins $A$ to $B$.
\end{enumerate}
\end{lemma}

\begin{proof}
Item~\ref{it:q-complementary-voltages}: Since $b\ne0$, the scalar master
identity \cref{lem:q-scalar-interface}\ref{it:q-master} gives
$h_u\ne0$ and $N-S_W-Np_u\ne0$ for every $u\in W$.
For distinct $u,v\in W$, subtract the identity indexed by $v$ from
that indexed by $u$:
\begin{align}
\label{eq:N-S}
0=(N-S_W-Np_u)h_u-(N-S_W-Np_v)h_v
 =(N-S_W)(h_u-h_v)-N(p_uh_u-p_vh_v).
\end{align}
If $h_u=h_v$, \eqref{eq:N-S} reduces to $N(p_u-p_v)h_u=0$, hence $p_u=p_v$.
On the other hand, if $p_u=p_v=t$, \eqref{eq:N-S} reduces to
$(N-S_W-Nt)(h_u-h_v)=0$; and since $N-S_W-Nt\ne 0$, we have $h_u=h_v$.
Therefore
\begin{align}
\label{eq:h-p-equiv}
h_u=h_v ~~\Longleftrightarrow~~ p_u=p_v.
\end{align}
When $p_u\ne p_v$, we have $h_u\ne h_v$, and the consistency identity
\cref{lem:q-scalar-interface}\ref{it:q-consistency} implies
$p_u+p_v=1$.  Thus at most two complementary voltage values occur.
They cannot all coincide: otherwise $h$ would be constant by the
equivalence above, and the zero sum identity \cref{lem:q-scalar-interface}\ref{it:q-zero-sum} would give $h=0$, contrary to
\cref{lem:q-scalar-interface}\ref{it:q-master} with $b\ne0$.
Thus exactly two values occur, denoted by $a$ and $1-a$ with
$0\le a<\frac12$.

Item~\ref{it:q-complementary-normal-form}: By the equivalence \eqref{eq:h-p-equiv},
$h$ is constant on each voltage class $A,B$; write its values as
$h_A,h_B$.  
\cref{lem:q-scalar-interface}\ref{it:q-zero-sum} gives
$|A|h_A+|B|h_B=0$.
Setting $\Delta:=h_A\ne0$ yields
$h_B=-\frac{|A|}{|B|}\Delta$.
Apply \cref{lem:q-scalar-interface}\ref{it:q-master} at $u\in A$, where $p_u=a$ and
$h_u=\Delta$.  Using $S_W=|A|a+|B|(1-a)$ and $N=|A|+|B|$ gives
\[
  (N+1)b=(S_W+Na-N)\Delta=|A|(2a-1)\Delta.
\]
Hence $b=R_0$ in \eqref{eq:complementary-pair-table}.
Substituting $b=R_0$, $h_A=\Delta$, and
$h_B=-\Delta\frac{|A|}{|B|}$ into
\cref{lem:q-scalar-interface}\ref{it:q-interpolation} reproduces the
$xA$, $yA$, $xB$, and $yB$ rows of
\eqref{eq:complementary-pair-table}.  Its formula for $F(uv)$ likewise
gives the $AB$, $AA$, and $BB$ rows; the remaining $xy$ entry equals
$b=R_0$.  Thus every vector for this fixed voltage partition is
determined by the single scalar $\Delta$.

Item~\ref{it:q-complementary-separation}: The coordinate calculation leading to
\eqref{eq:F-dist-AB} uses only the pair-value table
\eqref{eq:complementary-pair-table} and \eqref{eq:delta-coordinate}, not
the no-cross-edge assumption invoked later in the construction proof.
By item~\ref{it:q-complementary-normal-form}, the present $F$ has those
same pair values.  Thus, if $0<a<\frac12$, the formula gives, for
$u\in A$ and $v\in B$,
\[
\delta_F(u,v)^2-(p_u-p_v)^2\delta_F(x,y)^2
=\binom n2^{-1}\frac{2a(1-a)N^2}{|B|^2}\Delta^2>0.
\]
A positive $A$--$B$ edge would instead require equality by
\eqref{eq:edge-equality}, which is impossible.
If $a=0$, harmonicity at $u\in A$, where $p_u=0$, gives
$\sum_{z\in V}c_{uz}p_z=0$.  Every summand is nonnegative, and a
positive edge from $u$ to $v\in B$, where $p_v=1$, would make one
strictly positive. So there is no positive $A$--$B$ edge in this case either.
\end{proof}

\subsection{Existence criteria and kernel dimensions}

\begin{theorem}[Quadratic kernel classification]\label{thm:quadratic-classification}
Let $W_1,\dots,W_q$ be the connected components of $G[W]$.
Then
\begin{equation}\label{eq:quadratic-dimension}
  \dim\ker\!\left(\cK_{xy}^{(2)}\big|_{\cH_2}\right)
  =
  \begin{cases}
    q-1,
      &\text{if every $W_i$ is terminal-balanced and every $\rho_i=\frac12$ {\bf (the midpoint case)}},\\[1mm]
    1,
      &\text{if every $W_i$ is terminal-balanced and}
       \ \{\rho_1,\dots,\rho_q\}=\{a,1-a\}\\
      &\hspace{20mm}\text{for some $0\le a<\frac12$ {\bf (the complementary case)}},\\[1mm]
    0,
      &\text{otherwise.}
  \end{cases}
\end{equation}
\end{theorem}

\begin{proof}
\emph{Necessity and dimension upper bounds.}
Suppose $0\ne F\in\ker(\cK_{xy}^{(2)}|_{\cH_2})$, and let $p$ be its harmonic
voltage.  The alternatives $b=0$ and $b\ne0$ are covered by
\cref{lem:q-midpoint-alternative,lem:q-complementary-alternative}.

If $b=0$, \cref{lem:q-midpoint-alternative} gives $p\equiv\frac12$ on
$W$, so harmonicity implies $c_{xu}=c_{yu}$ at each $u\in W_i$.
Every $W_i$ is therefore terminal-balanced at $\frac12$.
Writing the componentwise constants as $\alpha_i$, its zero-sum condition is
$\sum_{i=1}^q |W_i|\alpha_i=0$.  Thus these vectors lie in a parameter space
of dimension $q-1$, proving the upper bound.

If $b\ne0$, \cref{lem:q-complementary-alternative} gives two distinct
voltages $a$ and $1-a$, both attained, with $0\le a<\frac12$, and no positive
edge between the corresponding classes.  Each connected component $W_i$
therefore lies in one class.  Writing its voltage as $\rho_i$, harmonicity
at $u\in W_i$ reduces to
$\rho_i c_{xu}=(1-\rho_i)c_{yu}$; hence every component is
terminal-balanced, with $\{\rho_1,\ldots,\rho_q\}=\{a,1-a\}$.
The harmonic $x$--$y$ voltage is uniquely determined by the graph and the
terminals.  Thus $a$ and the voltage classes $A=\{u\in W:p_u=a\}$ and
$B=\{u\in W:p_u=1-a\}$ are fixed, independently of $F$.  By \cref{lem:q-complementary-alternative}\ref{it:q-complementary-normal-form}, every such $F$ is specified by one
scalar, giving a dimension upper bound of one.

\emph{Sufficiency and attainment of the dimension bounds.}
Suppose first that every component is terminal-balanced at $\frac12$.
Define $p(x)=0$, $p(y)=1$, and $p\equiv\frac12$ on $W$.
Equation \eqref{eq:terminal-balance} makes $p$ harmonic.
For any constants $\alpha_i$ with
$\sum_{i=1}^q|W_i|\alpha_i=0$, let $\alpha\equiv\alpha_i$ on $W_i$
and define $F$ by \eqref{eq:q-midpoint-form}.
Its down-map sums are
$D_2F(x)=\sum_{u\in W}\alpha(u)=0$,
$D_2F(y)=-\sum_{u\in W}\alpha(u)=0$, and
$D_2F(u)=\alpha(u)-\alpha(u)=0$ for $u\in W$;
thus $F\in\cH_2$.
With $\Lambda:=\binom n2^{-1}\sum_{u\in W}\alpha(u)^2$,
\eqref{eq:delta-coordinate} yields
\[
  \delta_F(x,y)^2=4\Lambda,
  \qquad
  \delta_F(x,u)^2=\delta_F(u,y)^2=\Lambda
  \quad(u\in W).
\]
Because $\alpha$ is constant on each $W_i$,
\eqref{eq:q-midpoint-W-distance} gives zero $F$-distance on every
positive edge inside $W$.  The displayed terminal distances therefore
verify \eqref{eq:edge-equality} on every positive edge; the $xy$ edge,
if present, is automatic.  For nonzero $\alpha$, these distances also
show that $\phi_F=p$, so \cref{prop:MPL-equality,cor:operator-equality} give
$\cK_{xy}^{(2)}F=0$; the zero vector is in the kernel as well.

Conversely, let $F$ be a nonzero equality vector in the midpoint case.
Since $p\equiv\frac12$ on $W$,
\cref{lem:q-complementary-alternative} excludes $b\ne0$.
Thus \cref{lem:q-midpoint-alternative} gives
\eqref{eq:q-midpoint-form}, with $\alpha|_{W_i}=\alpha_i$ and
$\sum_{i=1}^q|W_i|\alpha_i=0$.  The parameters are unique, since
$\alpha_i=F(xu)$ for any $u\in W_i$.  Together with the preceding
construction, which also includes the zero vector, this identifies the
kernel linearly with
$
  \left\{(\alpha_i)_{i=1}^q\in\RR^q:
    \sum_{i=1}^q|W_i|\alpha_i=0\right\}$, which has dimension $q-1$.

Finally, suppose the terminal-balance parameters take exactly the two
values $a$ and $1-a$, with $0\le a<\frac12$.
Let $A$ and $B$ (both nonempty) be the unions of components with parameters $a$ and
$1-a$, respectively.
The function taking values $0,a,1-a,1$ on, respectively, $x,A,B,y$ is harmonic by
\eqref{eq:terminal-balance}.  No positive edge joins $A$ to $B$,
since distinct components of $G[W]$ have no such edge.
Thus \cref{lem:complementary-quadratic-construction} supplies a nonzero
quadratic kernel vector.  This attains the dimension upper bound of $1$
established above.
\end{proof}

\begin{remark}[Terminal-balance obstruction]
By \cref{thm:quadratic-classification}, if a component of $G[W]$ is not
terminal-balanced at any voltage, the quadratic equality kernel is zero.
In particular, a nonzero quadratic equality vector forces terminal balance
on every nonterminal component.
\end{remark}

\section{Degree lowering via zero-distance classes}
\label{sec:degree-lowering}

Throughout this section \cref{ass:top-degree-equality} is in force.
We retain the shorthands $W,N,p_u$, the interface
data $h,b$ from \eqref{eq:interface-data}, and $L=N-k+2>0$.
By \cref{prop:two-voltage-rigidity}, at most two nonterminal voltage values
occur.  We distinguish whether each voltage level is a single zero-distance
class or one level splits into several such classes.  These cases lead,
respectively, to complementary and midpoint graph geometry, from which the
quadratic classification supplies a degree-two equality vector.

\subsection{Single-class voltage levels: the complementary case}

We first consider the case in which each voltage value is a single zero-distance class.
If there is only one nonterminal voltage, the single-class hypothesis gives
$\delta_f(u,v)=0$ for all $u,v\in W$.  By \eqref{eq:transposition-distance},
every transposition $\tau_{uv}$ within $W$ fixes $f$ and hence $h$.  
These transpositions generate $\mathfrak S_W$, which
acts transitively on $(k-1)$-subsets, so $h\equiv c$.  Then $D_{k-1}^W h=0$
gives $0=(N-k+2)c=Lc$.  Since $L>0$, this forces $h=0$, contradicting
$\delta_f(x,y)>0$.  Hence there are exactly two nonterminal voltage classes.

Let $A,B$ be these classes, with respective voltages $p_A\neq p_B$.
By \cref{cor:zero-b-nonzero-h}, $h\neq0$.
Since vertices within each class have zero $\delta_f$-distance, their transpositions
fix $f$.  Consequently $h$ and $b$ are invariant under $\mathfrak S_A\times\mathfrak S_B$.
For $0\le j\le k-1$, let $h_j$ be the value of $h$ on a $(k-1)$-set containing
$j$ vertices of $B$, whenever that orbit exists.  Similarly, let $b_j$ be the value of
$b$ on a $(k-2)$-set containing $j$ vertices of $B$.

\begin{lemma}[Minimum sizes of the voltage classes]\label{lem:class-size}
If $h\neq0$, then $|A|,|B|\ge k-1$.
\end{lemma}

\begin{proof}
The standing range $k\le\lfloor \frac{n}{2}\rfloor$ gives $N=|A|+|B|\ge 2(k-1)$.
Suppose $|A|<k-1$, and set $j_0=k-1-|A|\ge1$.  The feasible indices of $h_j$
are then $j_0\le j\le k-1$ (in particular, $|B|\ge k-1$).
Choose $B'\subseteq B$ with $|B'|=j_0-1$, and put $Q=A\sqcup B'$.
Now $|Q|=|A|+j_0-1=k-2$, so the down-map identity \cref{lem:interface-identities}\ref{it:interface-down} applies:
\[
  0=(D_{k-1}^W h)(Q)=\sum_{z\in W\setminus Q} h_{Qz} = (|B|-j_0+1)h_{j_0}.
\]
Since $|B|-j_0+1>0$, we have $h_{j_0}=0$.
For each $j_0\le j<k-1$, choose a $(k-2)$-subset $Q'\subseteq W$
containing $j$ vertices of $B$.  Such a subset exists, and applying the same down-map identity to $Q'$ gives
\[
 0=(D_{k-1}^W h)(Q')= (|A|-k+2+j)h_j+(|B|-j)h_{j+1}.
\]
The coefficient $|B|-j$ is positive, so induction on $j$ shows that every
feasible $h_j$ vanishes, contradicting $h\ne0$.  Therefore $|A|\ge k-1$.
Interchanging $A$ and $B$ proves $|B|\ge k-1$ as well.
\end{proof}

By \cref{lem:class-size}, all orbit values $h_j$ for $0\le j\le k-1$
and $b_j$ for $0\le j\le k-2$ are defined.
Applying $D_{k-1}^W h=0$ to an arbitrary $(k-2)$-subset with $j$ vertices
of $B$ now gives the full recurrence
\begin{align}
 \label{eq:h-recurrence}
 (|A|-k+2+j)h_j+(|B|-j)h_{j+1}=0,\qquad 0\le j\le k-2.
\end{align}
Similarly, $D_{k-2}^W b=0$ gives
\begin{align}
\label{eq:b-recurrence}
  (|A|-k+3+j)b_j+(|B|-j)b_{j+1}=0,
  \qquad 0\le j\le k-3. 
\end{align}
For $0\le j\le k-3$, \cref{lem:interface-identities}\ref{it:interface-symmetry} gives
\begin{equation}\label{eq:b-h-difference}
  b_{j+1}-b_j=(p_A-p_B)h_{j+1}.
\end{equation}

\begin{lemma}[Complementary voltages]\label{lem:complementary}
The two unequal nonterminal voltages $p_A,p_B$ satisfy $ p_A+p_B=1$.
\end{lemma}

\begin{proof}
At $j=0$, the recurrences \eqref{eq:h-recurrence} and \eqref{eq:b-recurrence} become
$(|A|-k+2)h_0+|B|h_1=0$ and
$(|A|-k+3)b_0+|B|b_1=0$.
Substituting $b_1=b_0+(p_A-p_B)h_1$ from
\eqref{eq:b-h-difference} into the latter, and using
$L=N-k+2$, gives
\begin{align}
  h_0=-\frac{|B|}{|A|-k+2}h_1,
  \qquad
  b_0=-\frac{|B|(p_A-p_B)}{L+1}h_1.
\label{eq:h0-b0}
\end{align}
Since $h\neq0$, the recurrence \eqref{eq:h-recurrence} implies $h_1\neq0$.

Apply \cref{lem:interface-identities}\ref{it:interface-master} to a set $T$ consisting
of $k-2$ vertices of $A$, and a distinguished vertex $u\in A\setminus T$.
Then $b_T$ and $h_{uT}$ in \cref{lem:interface-identities}\ref{it:interface-master} correspond respectively to $b_0$ and $h_0$ in the current setting.
Substituting \eqref{eq:h0-b0} gives
$
  \frac{|B|L}{|A|-k+2}(p_A+p_B-1)h_1=0$.
Since $|B|>0, L>0$, $h_1\neq 0$, and $|A|-k+2>0$ by \cref{lem:class-size}, deduce that $p_A+p_B=1$.
\end{proof}

Relabel the two classes if necessary, and write $p_A=a$ and $p_B=1-a$ with $0\le a<\frac12$.

\begin{lemma}[No positive edges between voltage classes]\label{lem:no-cross-edge}
If $0\le a<\frac12$, then no positive edge joins $A$ to $B$.
\end{lemma}

\begin{proof}
Suppose on the contrary that there exist $u\in A$ and $v\in B$ with $c_{uv}>0$.
First suppose $a=0$, so that $p_A=0$ and $p_B=1$. Then
\[
\sum_{w\in V} c_{uw} [p(u)-p(w)] = c_{uv} [p(u)-p(v)] + \sum_{w\in V\setminus \{v\}} c_{uw} [p(u)-p(w)] \le -c_{uv} <0 
\]
which contradicts the harmonicity \eqref{eq:harmonic-voltage-equation} of the voltage $p$ at $u$.

Now suppose $0<a<\frac12$. By the
definition of the voltage,
$\delta_f(x,u)=p_u\delta_f(x,y)$ and $\delta_f(x,v)=p_v\delta_f(x,y)$, while
\eqref{eq:edge-equality} gives
$\delta_f(u,v)=\delta_f(x,y)(p_v-p_u)$.  Hence
$
  \delta_f(x,v)=\delta_f(x,u)+\delta_f(u,v)
$.
Apply \cref{lem:transposition-metric}\ref{it:transposition-equality} to the triple $(x,u,v)$ and to the
coordinate indexed by $yR$, where
$R\in\binomset{W\setminus\{u,v\}}{k-2}$.  By \cref{lem:interface-identities}\ref{it:interface-terminal},
$b_R-f(yuR)=p_uh_{uR}$ and $b_R-f(yvR)=p_vh_{vR}$.
These are respectively the $\mathbf a_{f,u}^{xv}$- and
$\mathbf d_f^{xv}$-coordinates in the equality criterion.
Their proportionality factor is $p_u/p_v$, and $p_u>0$;
thus $p_uh_{uR}=(p_u/p_v)p_vh_{vR}$, giving $h_{uR}=h_{vR}$.
The transposition $\tau_{uv}$ swaps $uR$ and $vR$, and fixes the $(k-1)$-subsets
which contain both or neither of $u,v$; thus it fixes $h$.  Together with
$\mathfrak S_A\times\mathfrak S_B$, this cross-transposition generates
$\mathfrak S_W$, so $h$ is constant on $(k-1)$-subsets.  Since $D_{k-1}^W h=0$, that constant must
be zero,
contradicting $h\neq0$.
\end{proof}

\begin{proposition}[Complementary-case graph geometry]
\label{prop:complementary-degree-lowering}
If each nonterminal voltage value is a single $\delta_f$-class, then every
connected component of $G[W]$ is terminal-balanced, and its terminal-balance
parameter belongs to $\{a,1-a\}$ for some $0\le a<\frac12$, with both values
attained.
\end{proposition}

\begin{proof}
By \cref{lem:complementary,lem:no-cross-edge}, the nonempty voltage classes $A,B$
have respective values $a,1-a$ and no positive edge joins them.
Thus each component $W_i$ of $G[W]$ lies in one class, where
$p\equiv\rho_i\in\{a,1-a\}$.  Harmonicity at $u\in W_i$ gives
$\rho_i c_{xu}=(1-\rho_i)c_{yu}$, the terminal-balance relation.
\end{proof}

\subsection{Split voltage levels: the midpoint case}

We now suppose that two distinct zero-distance classes occur at the same nonterminal
voltage.  Choose $u,v\in W$ with $p_u=p_v=t$ and $\delta_f(u,v)>0$.  For
$R\in\binomset{W\setminus\{u,v\}}{k-2}$ define $J_R:=h_{uR}-h_{vR}$.

\begin{lemma}[Same-voltage difference tensor]\label{lem:same-voltage-difference}
The tensor $J$ satisfies:
\begin{enumerate}[label=\textup{(\roman*)},ref=(\roman*)]
\item\label{it:same-voltage-down}
$D_{k-2}^{W\setminus\{u,v\}}J=0$.

\item\label{it:same-voltage-distance}
The distance identity
$\displaystyle 
  \delta_f(u,v)^2
  =\binom nk^{-1}\bigl[t^2+(1-t)^2\bigr]\sum_{R\in \binomset{W\setminus\{u,v\}}{k-2}} J_R^2$.
In particular $J\neq0$.
\item\label{it:same-voltage-rigidity}
For every $R\in\binomset{W\setminus\{u,v\}}{k-2}$, $(1-2t)J_R=0$; hence $t=\frac12$.
\end{enumerate}
\end{lemma}

\begin{proof}
Item~\ref{it:same-voltage-down}: Let
$Q\in\binomset{W\setminus\{u,v\}}{k-3}$.  Applying
\cref{lem:interface-identities}\ref{it:interface-down} to $uQ$ and $vQ$ gives
$h_{uvQ}+\sum_{z\notin uvQ}h_{uzQ}=0$ and
$h_{uvQ}+\sum_{z\notin uvQ}h_{vzQ}=0$.  Subtracting the two equations proves the claim.

Item~\ref{it:same-voltage-distance}: Apply the coordinate formula
\eqref{eq:delta-coordinate} to the pair $u,v$, and separate the indexing
sets $T\in\binomset{V\setminus\{u,v\}}{k-1}$ according to their intersection
with $\{x,y\}$.  If $T$ contains neither terminal, then
\cref{lem:interface-identities}\ref{it:interface-no-terminal} gives
\[
  f(uT)-f(vT)
  =(p_v-p_u)\bigl(f(xT)-f(yT)\bigr)=0,
\]
because $p_u=p_v=t$.  If $T=xyQ$ contains both terminals, then
\cref{lem:interface-identities}\ref{it:interface-symmetry} gives
$f(uT)-f(vT)=b_{uQ}-b_{vQ}=0$ for the same reason.

The remaining sets are $T=xR$ or $T=yR$, with
$R\in\binomset{W\setminus\{u,v\}}{k-2}$.  By
\cref{lem:interface-identities}\ref{it:interface-terminal}, their differences are
\[
\begin{aligned}
  f(xuR)-f(xvR)&=(1-t)(h_{uR}-h_{vR})=(1-t)J_R,\\
  f(yuR)-f(yvR)&=-t(h_{uR}-h_{vR})=-tJ_R.
\end{aligned}
\]
Each $R$ occurs once in each of these two classes of $T$.  Substituting all four
classes of $T$ into \eqref{eq:delta-coordinate} proves the distance identity.
Since $\delta_f(u,v)>0$, it follows that $J\ne0$.

Item~\ref{it:same-voltage-rigidity}: Fix
$R\in\binomset{W\setminus\{u,v\}}{k-2}$.  Apply
\cref{lem:interface-identities}\ref{it:interface-no-terminal} to the same
$k$-set $uvR$, first with the distinguished vertex $u$ and then with $v$.
Using \cref{lem:interface-identities}\ref{it:interface-terminal} in each case,
we obtain
\[
\begin{aligned}
 f(uvR)&=(1-t)f(xvR)+t f(yvR)
        =b_R+(1-2t)h_{vR},\\
 f(uvR)&=(1-t)f(xuR)+t f(yuR)
        =b_R+(1-2t)h_{uR}.
\end{aligned}
\]
Subtracting the two equations gives $(1-2t)J_R=0$ for every $R$.
Since $J\ne0$ by item~\ref{it:same-voltage-distance}, deduce that $t=\frac12$.
\end{proof}

\begin{lemma}[No second voltage in the split-level case]\label{lem:no-second-voltage}
If a nonterminal voltage level contains more than one zero-distance class, then
$p_w=\frac12$ for every $w\in W$.
\end{lemma}

\begin{proof}
By \cref{prop:two-voltage-rigidity}, at most one other voltage value can occur.
Suppose a second value $s\neq\frac12$ occurs on exactly $q_s\ge1$ vertices.  For
$R\in\binomset{W\setminus\{u,v\}}{k-2}$, let $j(R)$ be the number of those $q_s$ vertices
contained in $R$.  Apply \cref{lem:interface-identities}\ref{it:interface-master}
to $T=R$ and a vertex $u\in W\setminus R$, and again to the
same $T$ and a different vertex $v\in W\setminus R$. 
This produces two separate equations.
To simplify further, we note that
$p_u=p_v=\frac12$, so
$
S_W=\frac{N}{2}+q_s(s-\frac12)
$,
$S_R=\frac{k-2}{2}+j(R)(s-\frac12)$, and
\[
L-S_W- Lp_u + S_R = L-S_W-Lp_v + S_R = \frac{L}{2} - S_W+S_R = (q_s- j(R)) \left(\frac{1}{2}-s\right),
\]
which is the
common coefficient of $h_{uR}$ and $h_{vR}$ in the two equations.
Subtracting these equations yields
\begin{equation}\label{eq:J-support}
  \bigl(q_s-j(R)\bigr)\left(\frac12-s\right)J_R=0.
\end{equation}
Thus $J_R\neq0$ only if $R$ contains all $q_s$ vertices at voltage $s$.

If $q_s>k-2$ this is already impossible.  Otherwise, choose $R$ with $J_R\neq0$, remove
one voltage-$s$ vertex $s_0\in R$, and put $Q=R\setminus\{s_0\}$.  
Apply \cref{lem:same-voltage-difference}\ref{it:same-voltage-down}
to the set $Q$:
\[
 (D_{k-2}^{W\setminus\{u,v\}}J)(Q)
 =\sum_{z\in (W\setminus\{u,v\})\setminus Q} J_{Qz}=0.
\]
Only the term $z=s_0$ contains all $q_s$ vertices at voltage $s$;
every other summand vanishes by \eqref{eq:J-support}.  Thus $J_R=0$,
a contradiction.
\end{proof}

\begin{proposition}[Midpoint-case graph geometry]
\label{prop:midpoint-degree-lowering}
If a nonterminal voltage value contains more than one $\delta_f$-class,
then $G[W]$ has at least two connected components, and every component is
terminal-balanced at $\frac12$.
\end{proposition}

\begin{proof}
Choose two distinct zero-distance classes in the split voltage level.  
\cref{lem:no-second-voltage} implies that $p_w=\frac12$ for every $w\in W$.
For distinct $u,v\in W$ with $c_{uv}>0$, \eqref{eq:edge-equality}
forces $\delta_f(u,v)=0$.  Hence a connected component of $G[W]$ is contained in a
single zero-distance class.
In conjunction with the hypothesis, this implies that $G[W]$ has $q$ connected components with $q\ge 2$.

Since $p\equiv\frac12$ on $W$, harmonicity at each $u\in W$ reduces to
$c_{xu}=c_{yu}$.  Thus every connected component is terminal-balanced at
$\frac12$, as claimed.
\end{proof}

\subsection{Reduction to degree two}

\begin{theorem}[Reduction of higher-degree equality to degree two]\label{thm:all-degree}
Let $G=(V,c)$ be a finite connected weighted graph, fix $x\neq y$, and let
$3\le k\le \lfloor \frac{n}{2}\rfloor$.  If there exists a nonzero $f\in\cH_k$ such that
$
  \cK_{xy}^{(k)}f=0$,
then there exists a nonzero $F\in\cH_2$ such that
$
  \cK_{xy}^{(2)}F=0$.
\end{theorem}

\begin{proof}
Let $f\in\cH_k\setminus\{0\}$ satisfy $\cK_{xy}^{(k)}f=0$.
\cref{prop:two-voltage-rigidity} ensures at most two nonterminal voltage
values.  If some voltage level contains more than one $\delta_f$-class,
\cref{prop:midpoint-degree-lowering} gives at least two connected components of
$G[W]$, all terminal-balanced at $\frac12$.  Otherwise,
\cref{prop:complementary-degree-lowering} gives terminal-balanced components
whose parameters are exactly $\{a,1-a\}$ for some $0\le a<\frac12$.
Thus in either case the graph satisfies a nontrivial-kernel alternative of
\cref{thm:quadratic-classification}.  That theorem supplies a nonzero
$F\in\cH_2$ satisfying $\cK_{xy}^{(2)}F=0$.
\end{proof}

\section{Transfer between particle sectors}
\label{sec:sector-transfer}

To transfer the top-degree classification to an arbitrary particle sector, we
record the up-map realization of Hoeffding degree and its compatibility with the
exclusion residual $\cK_{xy}^{(\ell)}$.

For $0\le r<n$, define the one-step up operator
\[
  (U_r g)(\xi)
  :=\sum_{z\in \xi}g(\xi\setminus\{z\}),
  \qquad \xi\in\binomset V{r+1}.
\]
For $0\le j< \ell$, write $U_{j\to\ell}:=U_{\ell-1}\cdots U_j$, and set $U_{\ell\to\ell}=I$.
For $j\le\min(\ell,n-\ell)$, the degree-$j$ Hoeffding layer in the
$\ell$-particle sector is $U_{j\to\ell}\cH_j$.

\begin{lemma}[Up-map intertwining]\label{lem:up-map}
Let $1\le j\le\min(\ell,n-\ell)$.  The map $U_{j\to\ell}$ is injective on $\cH_j$ and
intertwines the vertex-permutation actions on its domain and codomain.  Consequently,
\begin{equation}\label{eq:K-intertwining}
  \cK_{xy}^{(\ell)}U_{j\to\ell}
  =U_{j\to\ell}\cK_{xy}^{(j)}.
\end{equation}
Thus a vector in $U_{j\to\ell}\cH_j$ is an equality vector if and only if its unique
preimage in $\cH_j$ is an equality vector.
\end{lemma}

\begin{proof}
Permutation equivariance is immediate from the definition of $U_r$, and therefore
\eqref{eq:K-intertwining} holds.

To prove injectivity, first note that for $1\le r<n$,
\begin{equation}\label{eq:DU-commutator}
  D_{r+1}U_r-U_{r-1}D_r=(n-2r)I.
\end{equation}
Indeed, for an $r$-set $T$ and a function $g$ on $r$-sets,
\begin{align*}
  (D_{r+1}U_rg)(T)
  &=(n-r)g(T)
   +\sum_{w\in T}\sum_{z\notin T}g(T\setminus\{w\}\sqcup\{z\}),\\
  (U_{r-1}D_rg)(T)
  &=r g(T)
   +\sum_{w\in T}\sum_{z\notin T}g(T\setminus\{w\}\sqcup\{z\}).
\end{align*}
Subtracting the two equations gives \eqref{eq:DU-commutator}.

If $g\in\cH_j$, set $g_0=g$ and $g_s:=U_{j+s-1}\cdots U_jg$ for $1\le s\le \ell-j$. 
Applying \eqref{eq:DU-commutator} inductively  gives
$D_{j+s}g_s=s(n-2j-s+1)g_{s-1}$.
Note that the coefficient $n-2j-s+1$ is positive because
$\ell\le n-j$.  
Hence $g_s=0$ implies $g_{s-1}=0$, and inductively $g=0$.
In particular, $ U_{j\to \ell} g=g_{\ell-j} =0$ implies $g=0$, and thus $U_{j\to\ell}$ is injective on $\cH_j$.

If $U_{j\to\ell}g$ lies in the kernel of $\cK_{xy}^{(\ell)}$, then
\eqref{eq:K-intertwining} gives
$U_{j\to\ell}(\cK_{xy}^{(j)}g)=0$.
The down map is equivariant under vertex permutations, as noted in the beginning of \cref{sec:interface}.  
Since $\cK_{xy}^{(j)}$ is a linear combination of the
identity and transpositions, this yields
$
  D_j\cK_{xy}^{(j)}=\cK_{xy}^{(j-1)}D_j$.
As $D_jg=0$, we have $D_j\cK_{xy}^{(j)}g=0$, hence
$\cK_{xy}^{(j)}g\in\cH_j$.  The injectivity of $U_{j\to\ell}$ on
$\cH_j$ therefore gives $\cK_{xy}^{(j)}g=0$.
The reverse implication follows from \eqref{eq:K-intertwining}.
\end{proof}

\begin{lemma}[Down-map surjectivity]\label{lem:down-map-surjective}
Let $E\subseteq V$ be a finite set.  For $1\le r\le \frac{|E|}{2}$, the down map
$D_r^E:\RR^{\binomset Er}\to\RR^{\binomset E{r-1}}$ is surjective.
\end{lemma}

\begin{proof}
For $r=1$, $D_1^E$ is the summation map onto functions on the unique $0$-subset, so it is
surjective.  Assume $r\ge2$.  Equip the subset-function spaces with the unnormalized counting
inner products.  For $0\le s<|E|$, denote the one-step up map on $E$ by $U_s^E$, where
$(U_s^E g)(\xi)=\sum_{z\in\xi}g(\xi\setminus\{z\})$ for $\xi\in\binomset E{s+1}$.
Then $D_r^E=(U_{r-1}^E)^*$.  
The commutator calculation in \eqref{eq:DU-commutator}, with $n$
replaced by $|E|$, gives
\[
  D_r^E U_{r-1}^E-U_{r-2}^E D_{r-1}^E=(|E|-2r+2)\,\mathrm{Id}.
\]
If $U_{r-1}^E g=0$, taking the inner product with $g$ yields
\[
  0=\|D_{r-1}^E g\|_{\mathrm{count}}^2
    +(|E|-2r+2)\|g\|_{\mathrm{count}}^2.
\]
The coefficient $|E|-2r+2$ is positive for $r\le \frac{|E|}{2}$, so $g=0$.  Thus $U_{r-1}^E$ is injective,
whence its adjoint $D_r^E$ is surjective.
\end{proof}

\begin{lemma}[Hoeffding--Specht decomposition across particle sectors]
\label{lem:hoeffding-specht-decomposition}
For $0\le j\le\lfloor \frac{n}{2}\rfloor$, the top-degree space $\cH_j$ satisfies
\begin{equation}\label{eq:top-specht-identification}
  \cH_j\cong S^{(n-j,j)}
\end{equation}
as an $\mathfrak S_V$-representation.  Moreover, for $1\le\ell\le n-1$,
put $d:=\min(\ell,n-\ell)$.  Then the decomposition
\begin{equation}\label{eq:hoeffding-direct-sum}
  \RR^{\Omega_\ell}
  =\bigoplus_{j=0}^{d}U_{j\to\ell}\cH_j
\end{equation}
is orthogonal.
\end{lemma}

\begin{proof}
Recall \eqref{eq:perm-Young}.
For $j\ge1$, the down map $D_j$ is
$\mathfrak S_V$-equivariant and surjective by \cref{lem:down-map-surjective}.  The
target contains exactly the summands with $0\le i\le j-1$.  Since each summand $S^{(n-i,i)}$ occurs with
multiplicity one, surjectivity forces $D_j$ to be nonzero, hence an isomorphism, on
each summand shared by the source and the target.  
Therefore $\ker D_j$ is precisely the remaining summand $S^{(n-j,j)}$, proving \eqref{eq:top-specht-identification}. 
The case $j=0$ is the constant representation.

For $1\le j\le d$, \cref{lem:up-map} shows that
$U_{j\to\ell}\cH_j$ is a nonzero invariant copy of $S^{(n-j,j)}$; for $j=0$ it
is the constant line.  Young's rule for the $\ell$-subset module, together with the
complement identification $\Omega_\ell\cong\Omega_{n-\ell}$ when
$\ell>\frac{n}{2}$, shows that $\RR^{\Omega_\ell}$ contains exactly one copy of each
$S^{(n-j,j)}$, $0\le j\le d$, and no others.  Hence the displayed sum is exhaustive.
The permutation action is orthogonal, and inequivalent irreducible summands are
orthogonal, which proves \eqref{eq:hoeffding-direct-sum}.
\end{proof}

\begin{corollary}[Equality above degree one if and only if quadratic equality]
\label{cor:degree-two-reduction}
Fix $G,x,y$ and an $\ell$-particle sector with $d:=\min(\ell,n-\ell)\ge2$.
Then
\[
  \ker\cK_{xy}^{(\ell)}\cap
  \biggl(\bigoplus_{j=2}^{d}U_{j\to\ell}\cH_j\biggr)\ne\{0\}
  \quad\Longleftrightarrow\quad
  \ker\cK_{xy}^{(\ell)}\cap U_{2\to\ell}\cH_2\ne\{0\}.
\]
\end{corollary}

\begin{proof}
By \cref{lem:hoeffding-specht-decomposition,lem:up-map}, each
$U_{j\to\ell}\cH_j$ is preserved by $\cK_{xy}^{(\ell)}$.
Thus a nonzero vector in the intersection on the left has a nonzero
component $g\in\ker\cK_{xy}^{(\ell)}\cap U_{j\to\ell}\cH_j$ for
some $2\le j\le d$.  If $j=2$, the conclusion follows immediately.
If $j\ge3$, write $g=U_{j\to\ell}h$ for the unique $h\in\cH_j$.
By \cref{lem:up-map}, $h$ is a nonzero equality vector.
\cref{thm:all-degree} supplies a nonzero $F\in\cH_2$ with
$\cK_{xy}^{(2)}F=0$, and \cref{lem:up-map} shows that
$U_{2\to\ell}F$ is a nonzero vector in the right-hand intersection.
The reverse implication follows from the inclusion of the $j=2$
summand in the direct sum.
\end{proof}

\subsection{Graph criterion for equality beyond degree one}

\begin{theorem}[Equality beyond degree one]
\label{thm:equality-beyond-degree-one}
Fix $G,x,y$ and an $\ell$-particle sector with
$d:=\min(\ell,n-\ell)\ge2$.  Put $W:=V\setminus\{x,y\}$,
and let $W_1,\dots,W_q$ be the connected components of $G[W]$.
Then
\begin{align}
\label{eq:ker-neq-0}
  \ker\cK_{xy}^{(\ell)}\cap
  \biggl(\bigoplus_{j=2}^{d}U_{j\to\ell}\cH_j\biggr)\ne\{0\}
\end{align}
if and only if either
\begin{enumerate}[label=\textup{(\alph*)},wide]
\item $q\ge2$ and every $W_i$ is terminal-balanced at $\frac12$; or
\item every $W_i$ is terminal-balanced, and the terminal-balance parameters are exactly
      $\{a,1-a\}$ for some $0\le a<\frac12$, with both values occurring.
\end{enumerate}
Moreover, the quadratic component in this sector satisfies
\[
  \dim\bigl(\ker\cK_{xy}^{(\ell)}\cap U_{2\to\ell}\cH_2\bigr)
  =\dim\ker\!\left(\cK_{xy}^{(2)}\big|_{\cH_2}\right).
\]
The common dimension is given by \cref{thm:quadratic-classification}.
\end{theorem}

\begin{proof}
By \cref{lem:up-map}, the restriction of $U_{2\to\ell}$ to
$\cH_2$ is injective and intertwines the residuals.  Consequently,
$
  \ker\cK_{xy}^{(\ell)}\cap U_{2\to\ell}\cH_2
  =U_{2\to\ell}\ker\!\left(\cK_{xy}^{(2)}\big|_{\cH_2}\right)$.
By \cref{cor:degree-two-reduction}, the intersection in \eqref{eq:ker-neq-0} is nonzero exactly when
$\ker(\cK^{(2)}_{xy}\big|_{\cH_2})\ne\{0\}$.  The graph alternatives are precisely those in
\cref{thm:quadratic-classification}, which also supplies the asserted dimension.
\end{proof}

\subsection{Examples: Weighted 4-cycles and complete bipartite graphs}

The first examples with equality beyond degree one occur on four vertices.
The next proposition records the complete weighted $4$-cycle family with
opposite terminals.

\begin{proposition}[Weighted $4$-cycles with opposite terminals]\label{prop:weighted-four-cycle}

Let $G$ be the $4$-cycle $x-u-y-v-x$ with positive conductances
$c_{xu}=\alpha$, $c_{uy}=\beta$, $c_{yv}=\gamma$, and $c_{vx}=\delta$.
Then $G$ supports equality beyond degree one for the terminal pair $x,y$ if and only if
$\beta\gamma=\alpha\delta$, in which case the quadratic equality kernel is one-dimensional.
\end{proposition}

\begin{proof}
Here $W=\{u,v\}$ and $G[W]$ has the two singleton components $\{u\}$ and $\{v\}$.
Their harmonic voltages are $p(u)=\beta/(\alpha+\beta)$ and
$p(v)=\gamma/(\gamma+\delta)$.
By \cref{thm:equality-beyond-degree-one}, equality beyond degree one occurs exactly when these
two terminal-balance parameters are complementary: $p(u)+p(v)=1$; 
in particular this subsumes the midpoint case. 
This constraint is equivalent to
$\beta\gamma=\alpha\delta$.  Since there are exactly two components,
\eqref{eq:quadratic-dimension} gives a quadratic equality kernel of dimension $1$.
\end{proof}

The uniformly weighted $4$-cycle is just one possible member of the equality family.  
For a non-uniformly weighted example, consider $(\alpha,\beta,\gamma,\delta)=(2,1,2,1)$, which gives harmonic
voltages $p(u)=\frac13$ and $p(v)=\frac23$.  

By contrast, if the terminals of a $4$-cycle are adjacent, \emph{i.e.,} $x-y-u-v-x$, then $G[W]$ is connected, so neither alternative
in \cref{thm:equality-beyond-degree-one} can hold.  The same observation shows
that adding a positive edge between $u$ and $v$ to an equality $4$-cycle with opposite terminals destroys equality beyond degree one.

\begin{remark}[Connection to the second eigenspace of $\mathcal{L}^{\rm IP}_G$]
Among all weighted $4$-cycles, only the uniformly weighted cycle exhibits the
exceptional second-eigenspace phenomenon described in
\cite[Theorem~2.4]{BelotserkovskiyChen}.
\end{remark}

A useful larger family is the weighted complete bipartite graph $K_{2,q}$, where the
terminals $x,y$ lie in one side of the bipartition, and the nonterminal vertices $u_1,\ldots,u_q$ lie in the opposite side.  
Write $a_i:=c_{xu_i}>0$, $b_i:=c_{yu_i}>0$, and
$\rho_i:=b_i/(a_i+b_i)$.
Each $\{u_i\}$ is a connected component of $G[W]$, with terminal-balance parameter
$\rho_i$.  \cref{thm:equality-beyond-degree-one} gives the complete alternative:
equality beyond degree one occurs either when every $a_i=b_i$ and $q\ge2$, in which case the
quadratic kernel has dimension $q-1$; or when the set of parameters $\rho_i$ is exactly
$\{a,1-a\}$ for some $0\le a<\frac12$, with both values occurring, in which case the
quadratic kernel is one-dimensional.  
Of course, when $q=2$ we recover the weighted $4$-cycles in \cref{prop:weighted-four-cycle}.  

\section{Higher-degree kernels and sectorwise classification}
\label{sec:higher-kernels}

We now determine the equality kernel in each Hoeffding degree.  For
$1\le k\le\lfloor \frac{n}{2}\rfloor$ set
$\cE_k:=\ker\!\left(\cK_{xy}^{(k)}\big|_{\cH_k}\right)$.
The $k=1$ case is treated in \cref{prop:degree-one-kernel} below, while the $k=2$ case is settled by
\cref{thm:quadratic-classification}.  For $k\ge3$, if
$\cE_k\ne\{0\}$, apply \cref{thm:equality-beyond-degree-one} in the
$k$-particle sector: there $\ell=d=k$ and $U_{k\to k}=I$, so the graph
must satisfy one of the midpoint or complementary voltage alternatives.
By \cref{lem:interface-reconstruction,cor:reconstruction-equality},
their higher-degree kernels can be constructed from compatible interface
data, with equality reduced to the positive-edge conditions within $W$.
In the midpoint case, these data form an invariant down-map kernel, while in the complementary case they obey an explicit orbit recurrence.
The resulting kernels are assembled into the particle-sector and
interchange-process classifications.

For a subgroup $\Gamma\le\mathfrak S_W$, where $W=V\setminus\{x,y\}$,
and a linear space $\mathcal U$ of functions on subsets of $W$, write
\[
\mathcal U^\Gamma:=\{g\in\mathcal U:g(\gamma\xi)=g(\xi)
\text{ for every }\gamma\in\Gamma\text{ and every admissible }\xi\}
\]
for the subspace of $\Gamma$-invariant functions.

\begin{lemma}[Dimension of the invariant down-map kernel]\label{lem:invariant-top-dimension}
Let $\Gamma$ be a permutation group acting on $W$, and 
$\beta_s$ denote the number of $\Gamma$-orbits on $s$-subsets of $W$.  If
$1\le r\le \frac{N}{2}$, then
\begin{equation}\label{eq:invariant-top-dimension}
  \dim\bigl(\ker D_r^W\bigr)^\Gamma=\beta_r-\beta_{r-1}.
\end{equation}
\end{lemma}

\begin{proof}
The space of $\Gamma$-invariant functions on $s$-subsets of $W$ has dimension
$\beta_s$, since such a function is determined by its values on the
$\Gamma$-orbits.  By \cref{lem:down-map-surjective}, $D_r^W$ is surjective
for $r\le \frac{N}{2}$.
Let $v$ be $\Gamma$-invariant on $(r-1)$-subsets, and choose $g$ on
$r$-subsets with $D_r^Wg=v$.  Its $\Gamma$-average $\bar g(\xi):=\frac1{|\Gamma|}\sum_{\gamma\in\Gamma}g(\gamma^{-1}\xi)$ is $\Gamma$-invariant, and
permutation-equivariance gives
$
D_r^W\bar g=\frac1{|\Gamma|}\sum_{\gamma\in\Gamma}
 (D_r^Wg)\circ\gamma^{-1}=v.
$
Thus the down map $D_r^W$ is surjective between the invariant function spaces, and its
kernel there is $(\ker D_r^W)^\Gamma$. 
The rank-nullity theorem yields \eqref{eq:invariant-top-dimension}.
\end{proof}

\subsection{Midpoint case}

Let $W_1,\dots,W_q$ be the connected components of $G[W]$, and set
$\Gamma:=\mathfrak S_{W_1}\times\cdots\times\mathfrak S_{W_q}$.
Let $\beta_s$ count the $\Gamma$-orbits on $s$-subsets of $W$.

\begin{theorem}[Higher-degree midpoint kernel]\label{thm:midpoint-higher-kernel}
Assume the midpoint voltage configuration: the harmonic $x$--$y$
voltage $p$ equals $\frac12$ on $W$.
For $3\le k\le\lfloor \frac{n}{2}\rfloor$, the map
\begin{equation}\label{eq:midpoint-kernel-map}
\cH_{k-1}(W)^\Gamma  \ni h \mapsto f_h \in \mathbb{R}^{\binom{V}{k}}, \qquad
f_h(\zeta)
=\left\{
\begin{array}{ll}
\frac12h_\xi, & \zeta=x\xi, ~\xi\in \binom{W}{k-1},\\
-\frac12h_\xi, & \zeta=y\xi, ~\xi\in \binom{W}{k-1},\\
0, & \text{otherwise},
\end{array}
\right.
\end{equation}
defines a linear isomorphism
$\cH_{k-1}(W)^\Gamma\cong\cE_k$.  
Consequently,
\begin{equation}\label{eq:midpoint-higher-dimension}
  \dim\cE_k=\beta_{k-1}-\beta_{k-2}.
\end{equation}
\end{theorem}

\begin{proof}
Let $f\in\cE_k$ be nonzero, and let $h$ and $b$ be its interface functions from
\eqref{eq:interface-data}.  By the midpoint hypothesis, $p\equiv\frac12$ on $W$, so
$S_W=\frac N2$ and $S_T=\frac{k-2}{2}$.  With $L=N-k+2$, the
coefficient of $h_{uT}$ in the master identity
\cref{lem:interface-identities}\ref{it:interface-master} is
$\frac12(L-N+k-2)=0$. 
Therefore
$(L+1)b_T=0$ for every $T$, hence $b=0$.  Using
\cref{lem:interface-identities}\ref{it:interface-terminal}--\ref{it:interface-no-terminal} now gives exactly
\eqref{eq:midpoint-kernel-map}.

If $c_{uv}>0$ with $u,v\in W$, then \eqref{eq:edge-equality} and
$p_u=p_v=\frac12$ give $\delta_f(u,v)=0$, hence $\tau_{uv}f=f$.
Since $\tau_{uv}$ fixes $x$ and $y$, it also fixes the interface function:
\[
  h_{\tau_{uv}\xi}
  =f(x\,\tau_{uv}\xi)-f(y\,\tau_{uv}\xi)
  =f(x\xi)-f(y\xi)=h_\xi.
\]
The positive-edge transpositions in each connected component $W_i$
generate $\mathfrak S_{W_i}$, so $h$ is $\Gamma$-invariant.
The down-map identity $D_{k-1}^W h=0$ in
\cref{lem:interface-identities}\ref{it:interface-down} then gives
$h\in\cH_{k-1}(W)^\Gamma$.  Thus every $f\in\cE_k$ has the prescribed form
\eqref{eq:midpoint-kernel-map} for some $h\in\cH_{k-1}(W)^\Gamma$.

Conversely, take $h\in\cH_{k-1}(W)^\Gamma$ and set $b=0$, as
required by the midpoint master identity (see the first paragraph of this proof).
If $h=0$,
then $f_h=0$, and the assertion is immediate. 
Hence assume $h\ne0$.
With the harmonic voltage $p\equiv\frac12$ fixed by the
hypothesis, we verify the compatibility conditions of
\cref{lem:interface-reconstruction}.  
\cref{lem:interface-identities}\ref{it:interface-down}
holds because $h\in\cH_{k-1}(W)$ and $b=0$.  In
\cref{lem:interface-identities}\ref{it:interface-symmetry}, both sides vanish because $b=0$ and all nonterminal
voltages are equal.  In \cref{lem:interface-identities}\ref{it:interface-consistency}, the factor
$1-p_u-p_v$ vanishes.  
\cref{lem:interface-identities}\ref{it:interface-master} holds because its $h_{uT}$ coefficient vanishes,
as computed above.
Thus \cref{lem:interface-reconstruction} produces $f_h\in\cH_k$ with
precisely the values in \eqref{eq:midpoint-kernel-map}.  By
\cref{cor:reconstruction-equality}, it remains only to check positive
edges within $W$.  On such an edge $\{u,v\}$ inside $W_i$,
$\Gamma$-invariance of $h$ and the formula \eqref{eq:midpoint-kernel-map} for $f_h$
give $\tau_{uv}f_h=f_h$.  Hence $\delta_{f_h}(u,v)=0$, matching
$p_u=p_v$, and thus $f_h\in\cE_k$.
The map is therefore surjective; it is injective since $f_h=0$
implies $h=0$ by \eqref{eq:midpoint-kernel-map}.

The orbit count $\beta_s$ is the number of tuples $(j_1,\dots,j_q)$ with
$0\le j_i\le |W_i|$ and $\sum_i j_i=s$.  Since
$k-1\le \frac{N}{2}$, \cref{lem:invariant-top-dimension} gives
\eqref{eq:midpoint-higher-dimension}.
\end{proof}

\subsection{Complementary case}

\begin{theorem}[Higher-degree complementary kernel]
\label{thm:complementary-higher-kernel}
Assume the complementary voltage configuration: the harmonic $x$--$y$
voltage $p$ is constant on each connected component of $G[W]$ and takes
exactly two values $a$ and $1-a$ on $W$, for some $0\le a<\frac12$.
Let $A:=\{u\in W:p_u=a\}$ and $B:=\{u\in W:p_u=1-a\}$.
For $3\le k\le\lfloor \frac{n}{2}\rfloor$, we have
\[
  \dim\cE_k
  =
  \begin{cases}
    1,&k-1\le\min\{|A|,|B|\},\\
    0,&k-1>\min\{|A|,|B|\}.
  \end{cases}
\]
When $\dim \cE_k=1$, $\cE_k$ is spanned by the function $f$ reconstructed
from the $\mathfrak{S}_A\times \mathfrak{S}_B$-invariant interface data constructed in the proof.
\end{theorem}

\begin{proof}
Let $f\in\cE_k$ be nonzero.  If two vertices in $A$ had positive
$\delta_f$-distance, then one voltage level would contain more than one zero-distance class,
and \cref{lem:no-second-voltage} would force every nonterminal voltage to equal
$\frac12$, a contradiction.  Thus $\delta_f(u,v)=0$ for every $u,v\in A$, and similarly
for every $u,v\in B$.  Since
$\delta_f(u,v)^2=\frac12\|f-\tau_{uv}f\|_k^2$, every such
transposition fixes $f$.  Hence $f$, and its corresponding interface functions $h$ and $b$ defined in \eqref{eq:interface-data}, are
invariant under $\mathfrak S_A\times\mathfrak S_B$.

Applying \cref{lem:class-size} at degree $k$ gives $|A|,|B|\ge k-1$ for every
nonzero equality vector.  Consequently, if $k-1>\min\{|A|,|B|\}$, the kernel is trivial.
So assume from now on that $|A|,|B|\ge k-1$.
Write $h_j$ for the value of $h$ on a $(k-1)$-subset containing $j$ points
of $B$.  The down-map identity $D_{k-1}^W h=0$ from \cref{lem:interface-identities}\ref{it:interface-down} gives \eqref{eq:h-recurrence}.
Both coefficients in \eqref{eq:h-recurrence} are positive for $0\le j\le k-2$, so the orbit values of $h$
are determined by $h_0$.  Let $b_j$ be the corresponding orbit value of $b$.
For a $(k-2)$-subset $T$ with $j$ points of $B$, choose $u\in A\setminus T$.
The master identity \cref{lem:interface-identities}\ref{it:interface-master}
then determines $b_j$ uniquely from $h_j$, since $N-k+3>0$.
Since $p$ is fixed, the uniqueness in
\cref{lem:interface-reconstruction}\ref{it:reconstruction-well-defined}
shows that $h_0$ determines $f$; hence $\dim\cE_k\le1$.
Moreover, $f\ne0$ implies $h\ne0$ by
\cref{cor:zero-b-nonzero-h}, hence $h_0\ne0$.
Taking $T\subset A$ and $u\in A\setminus T$ in
\cref{lem:interface-identities}\ref{it:interface-master} gives
$(N-k+3)b_0=(2a-1)(|A|-k+2)h_0\ne0$.
Thus every nonzero equality vector has $b\ne0$.

It suffices to prove existence of a nonzero vector $f\in \cE_k$.  Put
$L:=N-k+2$, $ \mu:=|A|-k+2$, and $\eta:=2a-1$.
Choose $h_0\ne0$, and define $h_1,\dots,h_{k-1}$ recursively by
\begin{equation}\label{eq:complementary-h-recurrence}
  (\mu+j)h_j+(|B|-j)h_{j+1}=0,
  \qquad 0\le j\le k-2.
\end{equation}
Set $h_\xi=h_j$ whenever $\xi\in\binomset W{k-1}$ contains $j$ vertices
of $B$.  For $T\in\binomset W{k-2}$ containing $j$ vertices of $B$, set
\begin{equation}\label{eq:complementary-b-formula}
  b_T=b_j,
  \qquad b_j:=\frac{\eta(\mu+j)}{L+1}h_j.
\end{equation}
The recurrence \eqref{eq:complementary-h-recurrence} gives $D_{k-1}^W h=0$.  Moreover,
for $0\le j\le k-3$,
\[
 (\mu+j+1)b_j+(|B|-j)b_{j+1}
 =\frac{\eta(\mu+j+1)}{L+1}
   \bigl[(\mu+j)h_j+(|B|-j)h_{j+1}\bigr]=0,
\]
so $D_{k-2}^W b=0$. We thus verified both down-map identities in \cref{lem:interface-identities}\ref{it:interface-down}. 

Now we verify the remaining hypotheses of \cref{lem:interface-reconstruction}.  Replacing a
vertex by another vertex in $A$ (resp.\@ $B$) does not change $h$ or $b$.  If one vertex
lies in $A$ and the other in $B$, then
\eqref{eq:complementary-h-recurrence} and
\eqref{eq:complementary-b-formula} give
\[
  b_{j+1}-b_j
  =\frac{\eta}{L+1}
    \bigl[(\mu+j+1)h_{j+1}-(\mu+j)h_j\bigr]
  =\eta h_{j+1}.
\]
Thus $b_{j+1}-b_j=\eta h_{j+1}$ for $0\le j\le k-3$, exactly as required by
\cref{lem:interface-identities}\ref{it:interface-symmetry}.  The identity
\cref{lem:interface-identities}\ref{it:interface-consistency} holds for
$T$ and distinct $u,v\in W\setminus T$: if $u,v$ both belong to $A$ (resp.\@ $B$),
the two $h$-values agree; if one of $\{u,v\}$ belongs to $A$ and the other to $B$, their voltages
sum to one.
Finally, we verify the master identity
\cref{lem:interface-identities}\ref{it:interface-master}.
Let $T$ contain $j$ vertices of $B$.  Since
$S_T=(k-2-j)a+j(1-a)$ and $S_W=|A|a+|B|(1-a)$,
the coefficient $L-S_W-Lp_u+S_T$ of $h_{uT}$ equals
$-\eta(\mu+j)$ if $u\in A\setminus T$, and
$\eta(|B|-j)$ if $u\in B\setminus T$.
In the two cases $h_{uT}=h_j$ and $h_{uT}=h_{j+1}$, respectively,
so \cref{lem:interface-identities}\ref{it:interface-master} reads
\begin{align}
 \label{eq:b_j-1}  (L+1)b_j-\eta(\mu+j)h_j&=0,
    \qquad u\in A\setminus T,\\
  \label{eq:b_j-2} (L+1)b_j+\eta(|B|-j)h_{j+1}&=0,
    \qquad u\in B\setminus T.
\end{align}
Equation \eqref{eq:b_j-1} is equivalent to \eqref{eq:complementary-b-formula}; substituting
\eqref{eq:b_j-1} into \eqref{eq:b_j-2} yields
\[
  \eta\bigl[(\mu+j)h_j+(|B|-j)h_{j+1}\bigr]=0,
\]
which holds by \eqref{eq:complementary-h-recurrence}.
Thus all compatibility conditions have been verified, and
\cref{lem:interface-reconstruction} reconstructs $f\in\cH_k$
from $(h,b,p)$, with $p$ the harmonic voltage prescribed by
the hypothesis of this theorem.

Since $h_0\ne0$, \cref{cor:reconstruction-equality} gives $f\ne0$
and reduces the equality test to positive edges within $W$.
The reconstruction formulas show that $f$ is
$\mathfrak S_A\times\mathfrak S_B$-invariant, since $h,b,p$
are invariant under this group.  Every positive edge within
$A$ or within $B$ therefore has zero $f$-distance, matching
its zero voltage drop; no positive edge joins $A$ to $B$.
By \cref{cor:reconstruction-equality}, $f\in\cE_k$,
completing the existence proof.
\end{proof}

\subsection{Top-degree and particle-sector kernels}

\begin{theorem}[Top-degree kernel classification]
\label{thm:top-degree-kernel-classification}
Let $W_1,\dots,W_q$ be the connected components of $G[W]$.
For every $2\le k\le\lfloor \frac{n}{2}\rfloor$, exactly the following alternatives occur:
\begin{enumerate}[label=\textup{(\roman*)},wide]
\item \label{item:top-degree-kernel-midpoint} If every $W_i$ is terminal-balanced at $\frac12$, let
      $\Gamma:=\mathfrak S_{W_1}\times\cdots\times\mathfrak S_{W_q}$ and let
      $\beta_s$ count the $\Gamma$-orbits on $s$-subsets of $W$.
      Then
      \[
        \dim\cE_k=\beta_{k-1}-\beta_{k-2}.
      \]
      Moreover, $\cE_k$ is canonically identified with
      $\cH_{k-1}(W)^\Gamma$ by \eqref{eq:midpoint-kernel-map},
      interpreted also at $k=2$.
\item \label{item:top-degree-kernel-complementary} If every $W_i$ is terminal-balanced, and the terminal-balance parameters are exactly
      $\{a,1-a\}$ for some $0\le a<\frac12$, let $A,B$ be the corresponding voltage
      classes.  Then
      \[
        \dim\cE_k
        =
        \begin{cases}
          1,&k-1\le\min\{|A|,|B|\},\\
          0,&k-1>\min\{|A|,|B|\}.
        \end{cases}
      \]
In the case $k\ge3$ and $\dim \cE_k=1$, the data in \eqref{eq:complementary-h-recurrence}--\eqref{eq:complementary-b-formula} together with \cref{lem:interface-reconstruction,cor:reconstruction-equality} give a spanning vector for $\cE_k$.  
In the case $k=2$, the vector in \eqref{eq:complementary-pair-table} spans $\cE_k$.
\item In all other cases, $\cE_k=\{0\}$.
\end{enumerate}
\end{theorem}

\begin{proof}
For $k=2$, the dimension statement is \cref{thm:quadratic-classification}.  In the
midpoint case, $\beta_1=q$ and $\beta_0=1$, so
$\beta_1-\beta_0=q-1$.  The proof of \cref{thm:quadratic-classification}
parametrizes this kernel by weighted-zero-sum component constants $\alpha$;
setting $h=2\alpha$ identifies their parameter space with $\cH_1(W)^\Gamma$
through the piecewise formula \eqref{eq:midpoint-kernel-map}, which is also valid at $k=2$.
In the complementary case, both voltage
classes are nonempty, so $1\le\min\{|A|,|B|\}$, $\dim \cE_k=1$, and $\cE_k$ is spanned by the vector \eqref{eq:complementary-pair-table}.

For $k\ge3$, alternatives (i) and (ii) are
\cref{thm:midpoint-higher-kernel,thm:complementary-higher-kernel}.
If neither graph alternative holds, \cref{thm:equality-beyond-degree-one}
with $\ell=k$ gives
$\ker\cK_{xy}^{(k)}\cap\bigl(\bigoplus_{j=2}^{k}
U_{j\to k}\cH_j\bigr)=\{0\}$.
Since $U_{k\to k}=I$, this forces $\cE_k=\{0\}$, proving the alternative (iii).
\end{proof}

To describe every particle sector, it remains to identify the degree-one equality
space.  Zhu \cite{Zhu26} established a nontrivial degree-one equality vector (which he called a ``linear observable'').
The next result shows that this vector spans $\cE_1$.

\begin{proposition}[Degree-one kernel]\label{prop:degree-one-kernel}
Let $p$ be the harmonic $x$--$y$ voltage and set $\bar p:=n^{-1}\sum_{v\in V}p(v)$.
Then
\begin{equation}\label{eq:degree-one-kernel}
   \cE_1
  =\operatorname{span}\{p-\bar p\}.
\end{equation}
\end{proposition}

\begin{proof}
The function $p-\bar p$ lies in $\cH_1$.  Its transposition distances are a constant
multiple of $|p(u)-p(v)|$, so \cref{prop:MPL-equality} shows that it is an equality
vector.  Conversely, let $g\in\cH_1$ be a nonzero equality vector.
For $r=1$ and $v\in V\setminus\{x,y\}$, apply
\eqref{eq:terminal-coordinate-proportionality} to $g$: the
$\mathbf a_{g,v}^{xy}$-coordinate is $g(x)-g(v)$, while the
$\mathbf d_g^{xy}$-coordinate is $g(x)-g(y)$.  Hence
$g(v)=(1-p(v))g(x)+p(v)g(y)$ for every vertex $v$, terminal or nonterminal.
Thus $g$ is affine in $p$, and the zero-sum condition defining $\cH_1$ leaves exactly
the one-dimensional span in \eqref{eq:degree-one-kernel}.
\end{proof}

\begin{theorem}[Kernel in each particle sector]
\label{thm:particle-sector-kernel}
Fix $1\le\ell\le n-1$ and put $d=\min(\ell,n-\ell)$. 
For $1\le j\le \lfloor \frac{n}{2}\rfloor$, let
$\cE_j$ be the kernels described in \cref{thm:top-degree-kernel-classification,prop:degree-one-kernel}.  
Then we have the orthogonal decomposition
\begin{equation}\label{eq:complete-sector-decomposition}
  \ker\cK_{xy}^{(\ell)}
  =\operatorname{span}\{\one\}
     \ \oplus\
   \bigoplus_{j=1}^{d}U_{j\to\ell}\cE_j,
\end{equation}
where $\one$ denotes the constant-one function.  
Consequently, using the orbit counts $\beta_s$ in the midpoint case, resp.\@
the voltage classes $A,B$ in the complementary case, as specified in
\cref{thm:top-degree-kernel-classification},
\begin{equation}\label{eq:sector-dimensions}
  \dim\ker\cK_{xy}^{(\ell)}
  =
  \begin{cases}
    1+\beta_{d-1},&\text{in the midpoint case},\\
    2+\min\{d-1,|A|,|B|\},&\text{in the complementary case},\\
    2,&\text{in all other cases}.
  \end{cases}
\end{equation}
\end{theorem}

\begin{proof}
By \cref{lem:up-map,lem:hoeffding-specht-decomposition}, the orthogonal
summands in \eqref{eq:hoeffding-direct-sum} are invariant under
$\cK_{xy}^{(\ell)}$.  For $1\le j\le d$, its kernel on
$U_{j\to\ell}\cH_j$ is exactly $U_{j\to\ell}\cE_j$.  In particular, the kernel intersection in
\cref{cor:degree-two-reduction} is the direct sum
$\bigoplus_{j=2}^{d}U_{j\to\ell}\cE_j$ when $d\ge2$.
Degree zero is the constant line, while \cref{prop:degree-one-kernel}
identifies the degree-one contribution.  This proves
\eqref{eq:complete-sector-decomposition}.

In the midpoint case, \cref{thm:top-degree-kernel-classification}\ref{item:top-degree-kernel-midpoint} and telescoping give the dimension
$
  2+\sum_{j=2}^{d}(\beta_{j-1}-\beta_{j-2})
  =2+\beta_{d-1}-\beta_0=1+\beta_{d-1}$.
In the complementary case, \cref{thm:top-degree-kernel-classification}\ref{item:top-degree-kernel-complementary} states that exactly
$\min\{d-1,|A|,|B|\}$ of the degrees $2,\dots,d$ contribute one dimension.  In all other cases no Hoeffding degree $j\ge2$ contributes.  This yields
\eqref{eq:sector-dimensions}.
\end{proof}

\begin{remark}
For $\ell\in \{0,n\}$, there is only one $\ell$-particle configuration.
Hence the residual $\cK^{(\ell)}_{xy}$ vanishes, and its kernel is the constant line.
\end{remark}

\subsection{Equality kernel for the interchange process}

The preceding \cref{thm:particle-sector-kernel} determines the equality kernel in
all exclusion particle sectors, including its component in each Hoeffding layer.
To pass to the interchange process, define the moving-particle residual on the
group algebra $\mathbb{R}[\mathfrak{S}_V]$ by
\[
  \cK_{xy}^{\mathrm{IP}}
  :=\cL_G^{\mathrm{IP}}-\ceff(x,y)(I-\tau_{xy}).
\]
On $S^\lambda$, the corresponding residual is $\cK_{G,xy}^\lambda$ from \eqref{eq:Q-L-K}.
Recall the orthogonal regular representation decomposition
\eqref{eq:regular-specht-decomposition}.
On the $\lambda$-isotypic summand $S^\lambda \otimes \mathbb{R}^{d_\lambda}$, $\cK_{xy}^{\mathrm{IP}}$ acts as
$\cK_{G,xy}^\lambda\otimes I$.

Set $\cE_0:=\RR$ for the constant representation.

\begin{theorem}[Complete equality kernel for the interchange process]
\label{thm:complete-interchange-kernel}
The interchange residual $\cK_{xy}^{\mathrm{IP}}$ defined above has the orthogonal decomposition
\begin{equation}\label{eq:complete-interchange-kernel}
  \ker\cK_{xy}^{\mathrm{IP}}
  \cong
  \bigoplus_{j=0}^{\lfloor \frac{n}{2}\rfloor}
    \cE_j\otimes\RR^{d_{(n-j,j)}}.
\end{equation}
In particular,
\begin{equation}\label{eq:interchange-kernel-dimension}
  \dim\ker\cK_{xy}^{\mathrm{IP}}
  =n+
   \sum_{j=2}^{\lfloor \frac{n}{2}\rfloor}
     \left(\binom{n}{j}-\binom{n}{j-1}\right)\dim\cE_j,
\end{equation}
where $\dim\cE_j$ for $j\ge2$ are given by
\cref{thm:top-degree-kernel-classification}.
\end{theorem}

\begin{proof}
For $1\le j\le\lfloor \frac{n}{2}\rfloor$, apply
\cref{thm:particle-sector-kernel} to the $j$-particle sector: its top
Hoeffding-layer kernel is $\cE_j$.  By
\eqref{eq:top-specht-identification}, this identifies
$\ker\cK_{G,xy}^{(n-j,j)}$ with $\cE_j$.  
By definition, $\cE_0=\RR$, and $\cE_1=\operatorname{span}\{p-\bar p\}$ by \cref{prop:degree-one-kernel}.
By \cref{thm:interchange-two-row-rigidity}, all irreducible summands with at least
three rows have trivial kernel.  Applying
\eqref{eq:regular-specht-decomposition} gives
\eqref{eq:complete-interchange-kernel}.
Since $\dim\cE_0=\dim\cE_1=1$, the $j=0,1$ summands contribute
$d_{(n)}+d_{(n-1,1)}=1+(n-1)=n$ dimensions.
For $j\ge2$, the hook-length formula gives
$d_{(n-j,j)}=\binom nj-\binom n{j-1}$, yielding
\eqref{eq:interchange-kernel-dimension}.
\end{proof}

\section*{Statement of AI use}

This author used ChatGPT 6--Sol to discover the two-row rigidity \cref{thm:interchange-two-row-rigidity}, the interface identities \cref{lem:interface-identities}, the voltage rigidity \cref{lem:k3-determinant,lem:four-point}, and the quadratic equality vector \eqref{eq:complementary-pair-table} in \cref{lem:complementary-quadratic-construction}.
In the writing process, multiple rounds of edit prompts were executed to ensure the correctness of the results, as well as their presentation.
The author manually checked the paper line-by-line and made the necessary edits before finalizing it. 
He takes full responsibility for the contents of this paper.

\end{document}